\documentclass[a4paper,leqno,11pt]{article}
\usepackage{latexsym}
\usepackage{amsmath}
\usepackage{amssymb}
\usepackage{enumerate}

\usepackage{theorem}

\theorembodyfont{\rmfamily}

\makeatletter

\@addtoreset{equation}{section}
\makeatother

\date{}
\title{On the existence of Generalized Unicorns on Surfaces
\footnote{
Mathematics Subject Classification (2000)\,:\,
 53B40, 
 Primary 53C60; 
 Secondary 53D35
 .}}

\author{S. V. Sabau, K. Shibuya and H. Shimada}
\begin{document}

\maketitle

\begin{abstract}
This paper addresses the problem of existence of generalized Landsberg structures on surfaces using the Cartan--K\"ahler Theorem and a 
Path Geometry approach.
\end{abstract}

\tableofcontents




\section{Introduction}
  \quad  A Finsler norm, or metric, on a real smooth, $n$-dimensional manifold
$M$ is a function $F:TM\to \left[0,\infty \right)$ that is positive and
smooth on $\widetilde{TM}=TM\backslash\{0\}$, has the {\it homogeneity property}
$F(x,\lambda v)=\lambda F(x,v)$, for all $\lambda > 0$ and all 
$v\in T_xM$, having also the {\it strong convexity} property that the
Hessian matrix
\begin{equation*}
g_{ij}=\frac{1}{2}\frac{\partial^2 F^2}{\partial y^i\partial y^j}
\end{equation*}
is positive definite at any point $u=(x^i,y^i)\in \widetilde{TM}$.\\
\quad The fundamental function $F$ of a Finsler structure $(M,F)$ determines and it is determined by the (tangent) {\it indicatrix}, or the total space of the unit tangent bundle of $F$
\begin{equation*}
\Sigma_F:=\{u\in TM:F(u)=1\}
\end{equation*}
which is a smooth hypersurface of $TM$.\\
\quad At each $x\in M$ we also have the {\it indicatrix at x}
\begin{equation*}
\Sigma_x:=\{v\in T_xM \ |\  F(x,v)=1\}=\Sigma_F\cap T_xM
\end{equation*}
which is a smooth, closed, strictly convex hypersurface in
$T_xM$. \\
   \quad A Finsler structure $(M,F)$ can be therefore regarded as smooth
hypersurface $\Sigma\subset TM$ for which the canonical projection
$\pi:\Sigma\to M$ is a surjective submersion and having the property
that for each $x\in M$, the $\pi$-fiber $\Sigma_x=\pi^{-1}(x)$ is
strictly convex including the origin $O_x\in T_xM$. We point out that the strong convexity condition of $F$ implies that the fiber $\Sigma_x$ is strictly convex, but the converse is not true (see \cite{BCS2000} for details on this point and a counterexample). \\
\quad A generalization of this notion is the {\it generalized Finsler structure} introduced 
by R. Bryant. In the two dimensional case a generalized Finsler structure is a 
coframing $\omega=(\omega^1,\omega^2,\omega^3)$ on a three dimensional manifold 
$\Sigma$ that satisfies some given structure equations (see \cite{Br1995}). 
By extension, one can study the generalized Finsler structure $(\Sigma,\omega)$ 
defined in this way ignoring even the existence of the underlying surface $M$. 
It was pointed out by C. Robles that in the case $n>2$, there will be no such 
globally defined coframing on the $2n-1$-dimensional manifold $\Sigma$. The reason 
is that even though the orthonormal frame bundle $\mathcal{F}$ over $M$ does admit 
a global coframing, it is a peculiarity of the $n=2$ dimensional case that  
$\mathcal{F}$ 
can be identified with $\Sigma$ (see also \cite{BCS2000}, p. 92-93 for concrete 
computations).  \\
\quad A generalized Finsler structure is {\it amenable} if the space of leaves $M$ 
of the foliation $\{\omega^1=0,\omega^2=0\}$ is 
differentiable manifold such that the canonical projection $\pi:\Sigma\to M$ is 
a smooth submersion. \\
\quad In order to study the differential geometry of the Finsler structure
$(M,F)$, one needs to construct the pull-back bundle
$(\pi^*TM,\pi,\Sigma)$ with the $\pi$-fibers $\pi^{-1}(u)$ diffeomorphic
to $T_xM$, where $u=(x,v)\in \Sigma$ (see
\cite{BCS2000}). In general this is not a principal bundle. \\
\quad By defining an orthonormal moving coframing on $\pi^*TM$ with respect to the
Riemannian metric on $\Sigma$ induced by the Finslerian metric $F$, the
moving equations on this frame lead to the so-called Chern
connection. This is an almost metric compatible, torsion free connection of the
vector bundle $(\pi^*TM,\pi,\Sigma)$. \\
\quad The canonical parallel transport $\Phi_t:T_xM\setminus 0 \to
T_{\sigma(t)}M\setminus 0 $, 
defined by the Chern connection along a
curve $\sigma$ on $M$, is a diffeomorphism that preserves the Finslerian
length of vectors. Unlike the parallel transport on a Riemannian
manifold,  $\Phi_t$ is not a linear isometry in general. \\
\quad This unexpected fact leads to some classes of special Finsler metrics. A Finsler metric whose parallel transport is a 
linear isometry is called a 
{\bf Berwald metric}, and one whose parallel transport is only a Riemannian  isometry is called a 
{\bf Landsberg metric} (see \cite{B2007} for a very good exposition). \\
\quad Equivalently, a Berwald metric is a Finsler metric whose Chern connection coincides with the Levi Civita 
connection of a certain Riemannian metric on $M$, in other words it is``Riemannian-metrizable". These are the closest 
Finslerian metric to the Riemannian ones. The connection is Riemannian, while the metric is not. However, in the two 
dimensional case, any Berwald structure is Riemannian or flat locally Minkowski, i.e. there are no geometrically 
interesting Berwald surfaces. \\
\quad Landsberg structures have the property that the Riemannian volume of the Finslerian unit ball is a constant. This 
remarkable property leads to a proof of Gauss-Bonnet theorem on surfaces \cite{BCS2000} and other interesting results.\
\
\quad Obviously, any Berwald structure is a Landsberg one. However, there are no examples of global Landsberg 
structures that are not Berwald. This is one of the main open problems in modern Finsler geometry.

{\bf Problem.} {\it Do there exist Landsberg structures that are not Berwald?}

\quad  The long time search for this kind of metric structures with beautiful properties, which everybody wanted to see 
but no one could actually get, makes D. Bao to call these metrics ``unicorns". \\
\quad On the other hand, on several occasions since 2002, R. Bryant claimed that there is plenty of {\it generalized} 
Landsberg structures on manifolds that are not Berwald. Moreover, he said that there are a lot of such generalized 
metrics depending on two families of two variables (see \cite{B2007}, p. 46--47). \\
\quad Even though from the first prophecy on the existence of generalized unicorns several years already passed, as far 
as we know, there is no proof or paper to confirm and develop further Bryant's affirmations. \\
\quad The purpose of this paper is two folded. First, we give a proof of the existence of generalized Landsberg structures 
on surfaces, which are not generalized Berwald structures and discuss their local amenability. \\
\quad Namely, we prove the following 

\quad {\bf Corollary 4.3.}\\
\quad {\it There exist non-trivial generalized Landsberg structures on a 3-manifold $\Sigma$.}

\quad Secondly, using a path geometry approach  we construct locally a generalized Landsberg structure by means of a 
Riemannian metric $g$ on the manifold of $N$-parallels $\Lambda$ (see 
\cite{Br1995} for a similar study of existence of generalized Finsler structures with $K=1$).
In the case when such Riemannian metric has its Levi-Civita connection $\nabla^g$ in a Zoll projective class 
$[\nabla]$ on $S^2$ it follows this generalized unicorn is in fact a classical one. We conjecture that this is always 
possible.\\
\quad In this way, even though we haven't explicitly computed yet the fundamental function $F:TM\to [0,\infty)$ of this 
Landsberg metric, our study gives 
an affirmative answer to the Problem posed above in the two dimensional case (see also \cite{Sz2008a}, \cite{M2008}, 
\cite{Sz2008b} for discussions on the existence of smooth unicorns in arbitrary dimension). Of course a proof for our 
conjecture in Section 9 remains to be given.

\quad Our method is based on an {\it upstairs} - {\it downstairs} gymnastics by moving between the base manifold and the 
total space of a fiber bundle. \\
\quad We give here the outline of our method in order to help the reader finding his way through the paper.\\
\quad We start by assuming  the existence of a generalized Landsberg structure $\{\omega^1,\omega^2,\omega^3\}$ on 
a 3-manifold $\Sigma$ and we perform first a coframe change (\ref{coframe_change}) by means of a function $m$ on $
\Sigma$ such that the new coframe $\{\theta^1,\theta^2,\theta^3\}$ has the properties:
\begin{enumerate}
\item it satisfies the structure equations (\ref{k_struct_eq});
\item its ``geodesic foliation'' $\{\theta^1=0,\theta^3=0\}$ coincides with the ``indicatrix foliation'' 
$\{\omega^1=0,\omega^2=0\}$ of the generalized Landsberg structure $(\Sigma,\omega)$.
\end{enumerate}
\quad Assuming these two conditions for $(\Sigma,\theta)$ we obtain a set of differential conditions for the function $m$ 
in terms of its directional derivatives with respect to the coframe $\omega$ given in Proposition 6.1, or, equivalently, 
in Proposition 6.2 if we start conversely. \\
    \quad Based on these, one can remark the following:
\begin{enumerate}
\item the function $m$ is invariant along the leaves of the foliation $\{\omega^2=0,\omega^3=0\}$, therefore, if we assume that 
$(\Sigma,\omega)$ is normal amenable, i.e. the leaf space of 
$\{\omega^2=0,\omega^3=0\}$ is a 2-dimensional differentiable manifold $\Lambda$, and the quotient projection 
$\nu:\Sigma\to \Lambda$ is a smooth submersion, then $m$ actually lives ``downstairs'' on this manifold rather than 
``upstairs'' on $\Sigma$ as initially expected;
\item If we realize $\{\theta^1,\theta^2,\theta^3\}$ as the canonical coframe of a Riemannian metric $g$ on $\Lambda$, 
then the function $k$ in (\ref{k_struct_eq}) is the lift of the Gauss curvature of the Riemannian metric $g$, hence the 
name ``curvature condition'' for (\ref{condC_up}) is motivated;
\item since we have constructed from the beginning the coframe $\theta$ on $\Sigma$ such that its geodesic foliation will 
generate the indicatrix leafs on $\Sigma$, if we could choose a Riemannian metric $g$ ``downstairs'' on $\Lambda$ all of 
whose geodesics are embedded circles, then the amenability of $(\Sigma,\omega)$ would be guaranteed. It is known 
that this kind of Riemannian metric exists and they are usually called Zoll metrics (see \cite{B1978}, \cite{G1976}). A 
more general concept is the Zoll projective structure $[\nabla]$ on $\Lambda$ (see \S 3.2 as well as \cite{LM2002}). 
These are projective equivalence classes of torsion free affine connections on 
$\Lambda$ whose geodesics are embedded circles in $\Lambda$. Moreover, under some very reasonable conditions 
they are metrizable by Riemannian metrics whose Levi-Civita connections $\nabla^g$ belong to the
 initial Zoll projective structure $[\nabla]$.
\end{enumerate}
    \quad All these imply that if we start ``downstairs'' with a Riemannian metric $g=u^2[(dz^1)^2+(dz^2)^2]$ on 
$\Lambda$, for some isothermal coordinates $(z^1,z^2)\in \Lambda$, where $u$ is a smooth function on $\Lambda$,  
then we can construct the $g$-orthonormal oriented frame bundle $\mathcal F(\Lambda)$ with its canonical 
coframe, say $\{\alpha^1,\alpha^2,\alpha^3\}$.\\
    \quad On the other hand, we set up a second order PDE system on $\Lambda$ for the functions $u, \bar m$ such that 
the lift $\widetilde m=\nu^*(\bar m)$ satisfies the conditions of Proposition 6.2. 
The Cartan-K\"ahler theorem tells us that such pairs of functions $(u,\bar{m})$ exist and they depend on 4 functions of 1 variable (Proposition 8.1). 
Then, by the coframe changing 
(\ref{inverse_coframe}) we obtain a new coframe $\widetilde \omega$ on the 3-manifold $\Sigma:=\mathcal F(\Lambda)$ 
which will satisfy the structure equations (\ref{Lands_struct_eq}) of a 
generalized Landsberg structure. The isothermal coordinates $(z^1,z^2)$ on $\Lambda$ and a homogeneous coordinate 
in the fiber of $\nu:\mathcal F(\Lambda) \to \Lambda$ over a point $z\in \Lambda$ will give a local form 
(\ref{normal_form}).\\ 
\quad The following diagram shows our {\it upstairs}-{\it downstairs} gymnastics.
\begin{equation*}
\begin{matrix}
Upstairs \quad  & (\Sigma,\omega)& \xrightarrow{m}& \quad (\Sigma,\theta)&\quad \equiv \quad & \quad\quad (\mathcal 
F(\Lambda),\alpha)& \xrightarrow{\widetilde m}& \Sigma:=(\mathcal F (\Lambda),\widetilde \omega)\\
& & & s^*\downarrow & & \nu^*\uparrow & & \\
Downstairs &  &  &\quad (\Lambda,g) &\quad \equiv \quad &(\Lambda,\widetilde g) & & 
\end{matrix}
\end{equation*}
\quad We use extensively the Cartan-K\"ahler theory in this paper in order to study the existence of integral manifolds of 
linear Pfaffian systems associated to PDE's upstairs as well as downstairs. The nontriviality of our generalized Landsberg structures can be achieved by 
choosing appropriate initial conditions for the integral submanifolds.\\
\quad The theory of exterior differential systems is one of the strongest tools
to study geometric structures. E. Cartan and other mathematicians 
reformulated various type of 
geometric structures by the exterior differential systems' terminology. However, very few essentially new results were 
obtained except for the work of R. Bryant, and few others (see \cite{Br et al 1991}, \cite{IL2003} and the references in 
these two fundamental books).\\
\quad In the present paper, the Cartan-K\"ahler theorem is essentially used to find the new 
geometric structures, namely generalized Landsberg structures. This shows the  
 usefulness and applicability of the theory of exterior differential systems. \\
 \quad For the concrete computations regarding Cartan-K\"ahler Theorem we have used the MAPLE package Cartan found in the Jeanne Clelland's home page (http://math.colorado.edu/ $\tilde{ }$ jnc/). We have found it extremely useful for checking this kind of computations.
\begin{center}
*
\end{center}

\quad Here is the structure of our paper. After a short survey of some basic notions of Finsler surfaces and generalized 
Finsler structures on surfaces in Section 2, we construct the linear Pfaffian exterior differential system in Section 3 whose integral 
manifolds are the sought structures. \\
\quad Using it we prove a local existence theorem for generalized Landsberg structures on surfaces that are not of 
Berwald type using the Cartan-K\"{a}hler theorem for linear Pfaffian systems in Section 4. Firstly, we assume the existence of generalized
 Landsberg structures on surfaces and build a linear Pfaffian system whose integral manifolds consist of the scalar invariants $I$ and $K$ 
 of the generalized Landsberg structure considered. Then Cartan-K\"{a}hler theorem tells us also that this 
kind of generalized structures depend on two arbitrary functions of two variables on $\Sigma$ (\S 4.1, \S 4.2). This proves Bryant's affirmations. \\
\qquad However, this discussion holds good under the assumption that generalized Landsberg structures exist.  We will show here more, namely,
we will study the involutivity of a Pfaffian system on $\Sigma$ whose integral manifolds consist of the coframe $\omega$ satisfying 
the structure equations \eqref{Lands_struct_eq} together with the scalar invariants $I$, $K$ satisfying the Bianchi identities \eqref{Lands_Bianchi}. 
This Pfaffian system is not a linear one, so we needed to prolong, but finally, 
Cartan-K\"{a}hler theorem tells us that these structures depend on 3 functions of 3 variables (\S 4.3). The degree of freedom is in this case 
higher than before, including the findings in \S 4.1, \S 4.2 as partial results. \\
\quad  
We discuss the 
local amenability of these structures in Section 5. \\
\quad Since the Cartan-K\"ahler theory is not very popular amongst the Finsler geometers, we introduce the basic 
notions and results in an Appendix. For the same reason, at the first use of the Cartan-K\"{a}hler theorem for linear Pfaffian systems, 
we present the computations in detail. Later uses
of the theorem in \S 4.3 and \S 8.2 show only the main formulas leaving the heavy computations to be verified by the reader. \\
\quad In order to obtain an amenable Landsberg structure on a 3-manifold 
$\Sigma$ we have considered in Section 6 a special coframe changing on $\Sigma$ constructed such that the indicatrix 
foliation of the initial Landsberg structure to coincide to the geodesic foliation of the new constructed structure. Moreover, 
this new coframe is realizable as the canonical coframe on the orthogonal frame bundle of a Riemannian surface (Section 7). \\
\quad Keeping all these in mind,  by inverting the procedure in Section 7 we have constructed in Section 8 a generalized Landsberg 
structure, on the total space $\mathcal F(\Lambda)$ of the 
orthonormal frame bundle of a Riemannian surface $(\Lambda,g)$, in terms of a  smooth function $\bar m$ on  $
\Lambda$. The Landsberg structure is not a Berwald one if  $\bar m$ is not constant.\\
\quad Finally, in Section 9, we discuss a possible way to show the existence of classical two dimensional unicorns. This problem is 
equivalent to the problem of finding a Riemannian metric $g$ that metrizes a Zoll projective class on $S^2$ and 
satisfies in the same time the PDE system (\ref{condL_2}), (\ref{condC_2}). We conjecture that this is always possible. \\
\quad Then, by construction the geodesic foliation 
$\{\alpha^1=0,\alpha^3=0\}$ of $g$ will foliate the 3-manifold $\mathcal F(\Lambda)$ by circles and the geodesic leaf 
space, say $M$, of the geodesic foliation naturally becomes a differentiable manifold and the leaf quotient mapping $\pi:
\mathcal F(\Lambda)\to M$ becomes a smooth submersion. In other words, we obtain a double fibration 
(see \S 3.1 and \S 3.2).\\
\quad Therefore, by our procedure it follows that this generalized Landsberg structure is amenable and its fibers $\pi^{-1}
(x)$ are compact, where $\pi:\mathcal F(\Lambda) \to M$, $x\in M$. \\
    \quad A simple argument will show that this generalized Landsberg structure is actually a classical Landsberg structure 
on the 2-manifold $M$, provided our conjecture is true. 

\quad {\bf Acknowledgments.} The authors would like to express their gratitude to Vladimir Matveev who pointed out an 
error in a previous version of the paper. We 
also thank to
David Bao, Gheorghe Pitis and Colleen Robles
for many useful discussions. We are also indebted to Keizo Yamaguchi for his valuable advises. 
Finally, we are grateful to the referee who pointed out the importance of the amenability of the generalized Landsberg structure and for many other helpful suggestions.
\section{Riemann--Finsler surfaces}

 	\quad We are going to restrict ourselves for the rest of the paper to the  two dimensional case. To be more precise, 
our manifold $\Sigma$ will be always 3-dimensional, and the manifold $M$ will be 2-dimensional, in the case it exists.\\
 	\quad{\bf Definition 2.1.}\quad A 3-dimensional manifold $\Sigma$ endowed with a coframing $
\omega=(\omega^1,\omega^2,\omega^3)$ which satisfies the structure equations 
\begin{equation}\label{finsler_struct_eq}
\begin{split}
d\omega^1&=-I\omega^1\wedge\omega^3+\omega^2\wedge\omega^3\\
d\omega^2&=-\omega^1\wedge\omega^3\\
d\omega^3&=K\omega^1\wedge\omega^2-J\omega^1\wedge \omega^3
\end{split}
\end{equation}
will be therefore called a {\it generalized Finsler surface},
where $I$, $J$, $K$ are smooth functions on $\Sigma$, called the invariants of the generalized Finsler structure $
(\Sigma,\omega)$ (see \cite{Br1995} for details).\\
	\quad As long as we work only with generalized Finsler surfaces, it might be possible that this generalized structure 
is not realizable as a classical Finslerian structure on a surface $M$. This imposes the following definition \cite{Br1995}.\\
	\quad{\bf Definition 2.2.}\ A generalized Finsler surface $(M,\omega)$ is said to be {\it amenable} if the leaf space $
\mathcal{M}$ of the codimension 2 foliation defined by the equations $\omega^1=0$, $\omega^2=0$ is a smooth surface 
such that the natural projection $\pi:\Sigma\to \mathcal{M}$   is a smooth submersion.\\
	\quad As R. Bryant emphasizes in  \cite{Br1995} the difference between a classical Finsler structure and a 
generalized one is global in nature, in the sense that {\it every generalized Finsler surface structure is locally 
diffeomorphic to a classical Finsler surface structure. 
}\\
	\quad The following fundamental result can be also found in \cite{Br1995}\\
	\quad{\bf Theorem 2.1.} \quad {\it The necessary and sufficient condition for a generalized Finsler surface  $
(\Sigma,\omega)$ to be realizable as a classical Finsler structure on a surface are
\begin{enumerate} 
\item the leaves of the foliation $\{\omega^1=0,\ \omega^2=0\}$ are compact;
\item it is amenable, i.e. the space of leaves of the foliation 
$\{\omega^1=0,\ \omega^2=0\}$ is a differentiable manifold $M$;
\item the canonical immersion $\iota:\Sigma\to TM$, given by 
$\iota(u)=\pi_{*,u}(\hat{e}_2)$, is one-to-one on each $\pi$-fiber $\Sigma_x$,
\end{enumerate}
where we denote by $(\hat{e}_1,\hat{e}_2, \hat{e}_3)$ the dual frame of the coframing $(\omega^1,\omega^2,\omega^3)$.
}\\
	\quad In the same source it is pointed out that if for example the $\{\omega^1=0,\ \omega^2=0\}$ leaves are not 
compact, or even in the case they are, if they are ramified, or if the curves  $\Sigma_x$ winds around origin in $T_xM$, 
in any of these cases, the generalized Finsler surface structure is not realizable as a classical Finsler surface.\\
	\quad An illustrative example found in \cite{Br1995} is the case of an amenable generalized Finsler surface such 
that the invariant $I$ is constant, however $I$ is not zero. This kind of generalized structure is not realizable as a Finsler 
surface because $I\neq 0$ means that the leaves of the foliation $\{\omega^1=0, \ \omega^2=0\}$ are not compact. 
Indeed, in the case $I^2<4$, the  $\pi$-fibers $\Sigma_x$ are logarithmic spirals in $T_xM$.\\
	\quad Let us return to the general theory of generalized Finsler structures on surfaces. 
By taking the exterior derivative of the structure equations (\ref{finsler_struct_eq}) one obtains the {\it Bianchi equations 
of the Finsler structure}:
\begin{equation*}
 \begin{split}
  & J=I_2\\
  & K_3+KI+J_2=0,
 \end{split}
\end{equation*}
 where we denote by $I_i$ the directional derivatives with respect to the coframing $\omega$, i.e.
\begin{equation*}
df=f_1\omega^1+f_2\omega^2+f_3\omega^3,
\end{equation*}
for any smooth function $f$ on $\Sigma$.\\
	\quad Taking now one more exterior derivative of the last formula written above,  one obtains the Ricci identities 
with respect to the generalized Finsler structure
\begin{equation*}
\begin{split}
& f_{21}-f_{12}=-Kf_3\\
& f_{32}-f_{23}=-f_1\\
& f_{31}-f_{13}=If_1+f_2+Jf_3.
\end{split}
\end{equation*}
	\quad{\bf Remarks.}
\begin{enumerate}
\item  Remark first that the structure equations of a Riemannian surface are obtained from (\ref{finsler_struct_eq}) by 
putting $I=J=0$. 
\item Since $J=I_2$, one can easily see that the necessary and sufficient condition for a generalized Finsler structure to 
be non-Riemannian is $I\neq 0$.  
\end{enumerate}
	\quad {\bf Definition 2.3.}\quad {\it A generalized Landsberg structure} on $\Sigma$ is a generalized Finsler 
structure $(M,\omega)$ such that $J=0$, or equivalently, $I_2=0$.\\
	\quad Remark that such a generalized
structure is characterized by the structure equations
\begin{equation}\label{Lands_struct_eq}
\begin{split}
d\omega^1&=-I\omega^1\wedge\omega^3+\omega^2\wedge\omega^3\\
d\omega^2&=-\omega^1\wedge\omega^3\\
d\omega^3&=K\omega^1\wedge\omega^2,
\end{split}
\end{equation}
and Bianchi identities
\begin{equation}\label{Lands_Bianchi}
\begin{split}
dI & =I_1\omega^1 \qquad \qquad +I_3\omega^3 \\
dK & =K_1\omega^1+K_2\omega^2-KI\omega^3,
\end{split}
\end{equation}
where $\omega=(\omega^1$, $\omega^2$,  $\omega^3)$ is a coframing on a certain 3-dimensional manifold $\Sigma$, 
and $I$ and $K$ are smooth functions defined on $\Sigma$. We will see that we actually need more, so we assume that 
the functions   $I$ and $K$ are analytic on $\Sigma$.\\
	\quad It is also useful to have the Ricci identities \cite{BCS2000} for the invariants $I$ and $K$. Indeed, taking first 
into account that
\begin{equation*}
K_{31}=-I_1K-IK_1,\quad K_{32}=-IK_2,\quad K_{33}=K(I^2-I_3),
\end{equation*}
we obtain
\begin{align*}
 & I_{12}=KI_3,          & K_{21}&-K_{12}=IK^2 \\
 & I_{32}=-I_1,          & K_{23}&=K_1-IK_2 \\
 & I_{31}-I_{13}=II_1,   & K_{13}&=-(2K_1I+KI_1+K_2).
\end{align*}
	\quad We are interested in studying the existence of non-trivial generalized Landsberg structures on $\Sigma$, i.e. 
generalized Landsberg structures that are not of Berwald type. \\
	\quad Recall the following definition.

	\quad {\bf Definition 2.4.} \quad A {\it generalized Berwald structure} is a generalized Finsler structure characterized 
by the structure equations (\ref{Lands_struct_eq}), and
\begin{equation*}
dI \equiv 0  \quad \mod\quad \omega^3,
\end{equation*}
or, equivalently,
\begin{equation*}
I_1=I_2=0,\qquad I_3\neq 0.
\end{equation*}
	\quad The reason we called Berwald structures (generalized or not) on surfaces {\it trivial} is given in the following 
rigidity theorem.\\

	\quad{\bf Theorem 2.2. Rigidity theorem for Berwald surfaces} \cite{Sz1981}\\
\quad {\it Let $(M,F)$ be a connected Berwald surface for which the Finsler structure 
$F$ is smooth and strongly convex on all of $\widetilde{TM}$.
\begin{enumerate}
\item If $K= 0$, then $F$ is locally Minkowski everywhere.
\item If $K\neq  0$, then $F$ is Riemannian everywhere.
\end{enumerate}
}
	\quad In other words, the only possible Berwald structures are either the flat locally Minkowski ones, or the 
Riemannian ones. Therefore the term {\it non-trivial} in the present paper addresses a Landsberg structure that is not 
locally Minkowski, nor Riemannian. Both of these are well studied trivial examples of Landsberg surfaces. \\
	\quad {\bf Remark.}\\
 It is interesting to remark that $I_1=0$ is not the only condition that makes a Lansdsberg to become a Berwald one. 

\quad Indeed, using the structure and the Ricci equations one can easily see that if for a Landsberg structure on a 
surface at least one of the following relation is satisfied
\begin{equation*}
 I_3=0, \quad K_2=0,
\end{equation*}
then that structure must be a Berwald one. \\
	\quad Remark also that the condition 
\begin{equation*}
K_1=0
\end{equation*}
does not necessarily imply triviality. In fact, all the generalized Landsberg structures in this paper satisfy this condition.
\section{Path Geometries} 
\subsection{Path geometries of a generalized Landsberg structure} 
\quad Recall from \cite{Br1997} that a {\it (classical) path geometry} on a surface $M$ is a foliation $\mathcal P$ of the 
projective tangent bundle $\mathbb P (TM)$ by contact curves, each of which is transverse to the fibers of the canonical 
projection $\pi:\mathbb P (TM)\to M$.\\
\quad Namely, let $\gamma:(a,b)\to M$ be a smooth, immersed curve, and let us denote by
$\hat{\gamma}:(a,b)\to \mathbb{P}(TM)$ 
its canonical lift 
to the projective tangent bundle $\pi:\mathbb{P}(TM)\to M$. Then, the fact that the canonical projection $\pi$ is a submersion implies that, for each line $L\in \mathbb{P}(TM)$, the linear map
\begin{equation*}
\pi_{*,L}:T_L \mathbb{P}(TM)\to T_xM,
\end{equation*}
is surjective, where $\pi(L)=x\in M$. Therefore
\begin{equation*}
E_L:=\bigl(\pi_{*,L}\bigr)^{-1}(L)\subset T_L \mathbb{P}(TM)
\end{equation*}
is a 2-plane in $T_L \mathbb{P}(TM)$ that defines a contact distribution and therefore a contact structure on $\mathbb{P}(TM)$.\\
\quad A curve on $\mathbb{P}(TM)$ is called {\it contact curve} if it is tangent to the contact distribution $E$. Nevertheless, the canonical lift $\hat{\gamma}$ to $\mathbb{P}(TM)$ of a curve $\gamma$ on $M$ is a contact curve.\\
\quad A {\it local path geometry} on $M$ is a foliation $\mathcal P$ of an open subset $U\subset \mathbb P (TM)$  by 
contact curves, each of which is transverse to the fibers of $\pi:\mathbb P (TM)\to M$.\\
\quad In the case there is a surface $\Lambda$ and a submersion $l:\mathbb P (TM) \to \Lambda$ whose fibers are the 
leaves of $\mathcal P$, then the path geometry will be called {\it amenable}.\\
\quad More generally, a {\it generalized path geometry} on a 3-manifold $\Sigma$ is a pair of transverse codimension 2 
foliations $(\mathcal P, \mathcal Q)$ with the property that the (unique) 2-plane field $D$ that is tangent to both foliations 
defines a contact structure on $\Sigma$. \\
\quad In the case when there is a surface $\Lambda_{\mathcal P}$ and a submersion 
$l_{\mathcal P}:\mathbb P (TM) \to \Lambda_{\mathcal P}$ whose fibers are the 
leaves of the foliation $\mathcal P$, then the generalized path geometry 
$(\mathcal P, \mathcal Q)$ 
 will be called $\mathcal P$-{\it amenable}. A $\mathcal Q$-{\it amenable}   generalized path geometry 
$(\mathcal P, \mathcal Q)$ is defined in a similar way.\\
\quad One can easily see that a classical path geometry on $\Sigma=\mathbb P (TM)$ is a special case of generalized 
path geometry where the second foliation $ \mathcal Q$ is taken to be the fibers of the canonical projection $\pi:\mathbb 
P (TM)\to M$.\\
\quad In the case of a Landsberg structure on a 3-manifold $\Sigma$, we can define two kinds of generalized path 
geometries.\\
\quad We can consider
\begin{enumerate}
\item  $ \mathcal P := \{\omega^1=0,\omega^3=0\}$ the {\it ``geodesic'' foliation} of $\Sigma$, i.e. the leaves are curves on 
$\Sigma$ tangent to $\hat e_2$;
\item  $ \mathcal Q := \{\omega^1=0,\omega^2=0\}$ the {\it ``indicatrix'' foliation} of $\Sigma$, i.e. the leaves are curves on 
$\Sigma$ tangent to $\hat e_3$;
\item  $ \mathcal R := \{\omega^2=0,\omega^3=0\}$ the {\it ``normal'' foliation} of $\Sigma$, i.e. the leaves are curves on $
\Sigma$ tangent to $\hat e_1$.
\end{enumerate}
	\quad We can consider now the generalized path geometries 
\begin{equation*}
 \mathcal G_1=(\mathcal P, \mathcal Q),\qquad \mathcal G_2=(\mathcal R, \mathcal Q).
\end{equation*}
	\quad Remark that on the case of $\mathcal G_1$, the 2-plane field $D_1=<\hat e_2,\hat e_3>$ defines indeed a 
contact structure on $\Sigma$. To verify this we need a contact 1-form $\eta$ on $\Sigma$ such that $D_1=\ker \eta$. By 
definition it follows that $\eta$ has to be 
\begin{equation*}
 \eta=A\omega^1
\end{equation*}
and we have
\begin{equation*}
 \eta\wedge d\eta=A^2\omega^1\wedge\omega^2\wedge\omega^3.
\end{equation*}
	\quad Therefore $\eta$ is a contact form on $\Sigma$ if and only if $A\neq 0$, so  $\mathcal G_1$ is a well defined 
path geometry on $\Sigma$. \\
	\quad We can do the same discussion for $\mathcal G_2$, where the 2-plane field is $D_2=<\hat e_1,\hat e_3>$. 
As above, we look for
$\eta$ such that  $D_2=\ker \eta$, so we must have 
\begin{equation*}
 \eta=A\omega^2,
\end{equation*}
and a simple computation shows that again
\begin{equation*}
 \eta\wedge d\eta=A^2\omega^1\wedge\omega^2\wedge\omega^3.
\end{equation*}
	\quad Therefore, again $\eta$ is a contact form on $\Sigma$ if and only if $A\neq 0$, and again  $\mathcal G_2$ is 
a well defined path geometry on $\Sigma$. \\
	\quad Recall also from the same \cite{Br1997} that {\it every generalized path geometry is always identifiable with a 
local path geometry on a surface}. Indeed, for a $u\in \Sigma$, let $U\subset \Sigma$ be an open neighborhood of $u$ 
on which the foliation $\mathcal Q$ is locally amenable, i.e. there exist a smooth surface $M$ and a smooth surjective 
submersion $\pi:U\to M$ such that the fibers of $\pi$ are the leaves of $Q$ restricted to $U$. Remark that this is always 
possible (for example due to Frobenius theorem) and that $M$ and $\pi$ are uniquely determined by $U$ up to a 
diffeomorphism.\\
	\quad A natural smooth map $\nu:U\subset \Sigma\to \mathbb P (TM)$, 
which makes the following diagram commutative,  
\begin{equation*}
\begin{split}
 & \qquad \ \ \quad \nu \\
 & U\subset  \Sigma \longrightarrow \mathbb P (TM) \\
\pi & \downarrow \qquad \swarrow \\
 & M
\end{split}
\end{equation*}
can be defined as follows
\begin{equation*}
\nu(u)=\pi_{*,u}(T_u\mathcal P),
\end{equation*}
for any $u\in U$. This application is  well defined because $\pi_{*,u}(T_u\mathcal P)$ is a 
1-dimensional subspace of $T_{\pi(u)}M$, and therefore an element of $\mathbb P (T_{\pi(u)}M)$.  \\
	\quad For the generalized path geometry $\mathcal G_1=(\mathcal P, \mathcal Q)$ we put
\begin{equation*}
 \nu_1:U\subset \Sigma\to \mathbb P (TM),\qquad \nu_1(u)=\pi_{*,u}(\hat e_2),
\end{equation*}
and for the 
generalized path geometry $\mathcal G_2=(\mathcal R, \mathcal Q)$ we put
\begin{equation*}
 \nu_2:U\subset \Sigma\to \mathbb P (TM),\qquad \nu_2(u)=\pi_{*,u}(\hat e_1).
\end{equation*}
	\quad Remark that because the foliations $\mathcal P$, $\mathcal Q$ and $\mathcal R$ are all transverse to  each 
other, it follows again that $\pi_{*,u}(T_u \mathcal P)$ and $\pi_{*,u}(T_u \mathcal R)$ are 1-dimensional subspaces in 
$T_{\pi(u)}M$, i.e. $\nu_1$, $\nu_2$ are immersions and therefore local diffeomorphisms.

\subsection{Zoll projective structures}
	\quad A classical example of a path geometry on a 3-manifold is the path geometry of a Riemannian metrizable Zoll 
projective structure. This is not only an example of path geometry, but it will be very useful in the construction of a 
non-trivial Landsberg structure.\\
	\quad Recall  that a Riemannian metric $g$ on a smooth manifold $\Lambda$ is called a {\it Zoll metric} if all its 
geodesics are simple closed curves of equal length. See \cite{B1978} for basics of Zoll metrics as well as 
\cite{G1976} for the abundance of Zoll metrics on $S^2$.\\
	\quad We will use in the present paper a more general notion, namely Zoll projective structures. Our exposition 
follows closely \cite{LM2002}.\\
	\quad{\bf Definition 3.1.} If $\nabla$ is a torsion free affine connection on a smooth manifold $\Lambda$, then the 
projective class $[\nabla]$ of $\nabla$ is called a {\it Zoll projective structure} if the image of any maximal geodesic of $
\nabla$ is an embedded circle $S^1\subset \Lambda$.\\
	\quad Given a Zoll projective structure $[\nabla]$ on $\Lambda$, the canonical lift of its geodesics will provide the 
geodesic foliation $\mathcal P$ on the projectivized tangent  bundle $\mathbb P(T\Lambda)$ which foliates 
$\mathbb P(T\Lambda)$ by circles. Let $M$ be the leaf space of the geodesic foliation $\mathcal P$ of a Zoll projective 
structure.\\
	\quad It can be shown that any Zoll projective structure $[\nabla]$ on a compact orientable surface $\Lambda$
is {\it tame}, namely each leaf of its geodesic foliation on $\mathbb P(T\Lambda)$ has a neighborhood which is 
diffeomorphic to $\mathbb R^2\times S^1$, such that each leaf corresponds to a circle of the form
$\{u\}\times S^1$, for any $u\in \mathbb P(T\Lambda)$.\\
	\quad This implies further that the leaf space $M$ of a Zoll projective structure $[\nabla]$ on a compact orientable 
surface $\Lambda$ has a canonical structure of differentiable manifold such that the quotient map $\pi:\mathbb P(T
\Lambda) \to M$ becomes a submersion. We obtain therefore the following  {\it double fibration} of a Zoll projective 
structure.
\begin{equation*}
\begin{split}
 & \qquad \mathbb P (T\Lambda) \\
 & \nu  \swarrow \qquad \searrow \pi\\
 & \Lambda\qquad\qquad \quad M
\end{split}
\end{equation*}
	\quad Let us assume from now $\Lambda=S^2$. It is natural to ask when a given Zoll projective structure 
$[\nabla]$ on $S^2$ can be represented by the Levi-Civita connection of a Riemannian metric $g$ on $\Lambda=S^2$. \\
	\quad The answer is given in Theorem 4.8. of \cite{LM2002}, p. 512. We are not going to state or to prove this 
theorem here because it will take too much space to define all the notions that are involved. Instead, we are going to 
describe the construction of the Riemannian metric $g$ that represents a Zoll projective structure, in the case such a 
metric exists. It is clear from \cite{LM2002} that the set of Riemannian metrizable Zoll projective structures 
is not empty, so we can assume the existence of Riemannian metrizable Zoll projective structures $[\nabla]$ on $S^2$. \\
	\quad Let us consider the isothermal local coordinates $(z^1,z^2)$ on $S^2$ induced from the Zoll projective 
structure (the concrete construction can be found in \cite{LM2002}, p. 513), and let
%
%
\begin{equation*}
g=u^2\Bigl[ (dz^1)^2+(dz^2)^2 \Bigr],
\end{equation*}
be the metric on $S^2$ in these coordinates,
where $u$ is a smooth function. If one puts
\begin{equation*}
\gamma=d\log u,
\end{equation*}
then the Levi-Civita connection $\nabla^g$ of  the Riemannian metric $g$ belongs to the Zoll projective structure $
[\nabla]$ if
\begin{equation}\label{Gamma}
\Gamma_{kl}^j=\gamma_k\delta_l^j+\gamma_l\delta_k^j-\gamma^j\delta_{kl},
\end{equation}
where $\gamma=\gamma_1dz^1+\gamma_2dz^2$, and $\Gamma_{kl}^i$ are the Christoffel symbols of the 
Zoll projective structures $[\nabla]$, i.e. 
\begin{equation*}
\Gamma_{jk}^i=\Bigl< dz^i,\nabla_{\frac{\partial}{\partial z^j}}\frac{\partial}{\partial z^k}\Bigr>
\end{equation*}
for a connection $\nabla$ in the Zoll projective structure $[\nabla]$, and $\gamma^j=g^{ji}\gamma_i$.\\
	\quad It follows that for a given Zoll projective structure $[\nabla]$ we obtain
\begin{equation}\label{small_gamma}
\gamma_i=\frac{1}{2}\Bigl(\Gamma_{i1}^1+\Gamma_{i2}^2 \Bigr),\qquad i=1,2.
\end{equation}
 	\quad If we denote by $R$ the Gauss curvature of $g$, then taking into account that 
$\gamma_i=\frac{1}{u}\frac{\partial u}{\partial z^i}$, it follows
\begin{equation}
R=-\frac{1}{u^2}\textrm{div} \gamma, 
\end{equation}
where we put $\textrm{div} \gamma=\frac{\partial\gamma_1}{\partial z^1}+\frac{\partial\gamma_2}{\partial z^2}$.\\
	\quad If we denote by $\{\alpha^1,\alpha^2,\alpha^3\}$ the canonical coframe on the bundle of $g$-orthonormal 
frames on $\Lambda$ then
$\mathcal G=(\mathcal P, \mathcal Q)$
is a path geometry on $\mathbb P(T\Lambda)$, where $\mathcal P:=\{\alpha^1=0,\alpha^2=0\}$ is the
geodesic foliation of $g$ and  $\mathcal Q:=\{\alpha^1=0,\alpha^3=0\}$. 
\section{The Cartan--K\"ahler theory}

\subsection {A linear Pfaffian system on generalized Landsberg surfaces}
\quad This section and the following one are motivated by Bryant's prophecy on the existence of generalized unicorns that we mentioned already in Introduction. Since the Finsler geometry community is familiarized with his statements, we will give here our interpretation of it. We point out however, that the discussion following hereafter does  not imply the existence of non-trivial generalized unicorns.  This will be shown only in section 4.3 in a different setting.\\
	\quad In order to make use of the  
	Cartan-K\"{a}hler theory, we are going to construct an exterior differential system associated to the coframe 
	$(\omega^1,\omega^2,\omega^3)$ that satisfies (\ref{Lands_struct_eq}), 
(\ref{Lands_Bianchi}).\\ 
\quad In this section we {\it assume} the existence of three linear independent one forms 
$(\omega^1,\omega^2,\omega^3)$ on the 9-dimensional manifold 
$\widetilde{\Sigma}=\Sigma\times \mathbb R^2\times \mathbb R^4$ that satisfy the structure equations (\ref{Lands_struct_eq}), where we consider the free coordinates $(I,K)\in \mathbb R^2$, and 
$(I_1,I_3,K_1,K_2)\in \mathbb R^4$, 
 and study the degree of freedom of the scalar functions $I$ and $K$.\\
	\quad First, we consider the following 1-forms
\begin{equation}\label{Lands_Pfaffian}
\begin{split}
\theta^1: & = dI-I_1\omega^1-I_3\omega^3 \\
\theta^2: & = dK -K_1\omega^1-K_2\omega^2+KI\omega^3,
\end{split}
\end{equation}
and let us denote by $\mathcal{I}$ the differential ideal generated by $\{\theta^1,\theta^2\}$. We also denote
\begin{equation*}
\begin{split}
\Omega & :=\omega^1\wedge\omega^2\wedge\omega^3, \\
J & :=\{\theta^a,\omega^i\}, \\
I& :=\{\theta^a\},
\end{split}
\end{equation*}
where {\it a}=1,2, {\it i}=1,2,3.\\
	\quad We will use the same letter $I$ for the invariant of a (generalized)
Finsler structure as well as for the set of 1-forms
$\theta^1,\theta^2$. We hope that this will not lead to any confusion.\\
	\quad In order to use the Cartan--K\"{a}hler theory we are going to consider the pair $(I,J)$ as an EDS with 
independence condition on a certain manifold $\widetilde{\Sigma}$ to be determined later.  We consider $dI$ and $dK$ 
as linearly independent 1-forms on the manifold $\widetilde{\Sigma}$.\\
	\quad By exterior differentiation  of $\{\theta^1,\theta^2\}$ we obtain
\begin{equation*}
\begin{split}
d\theta^1 & =-dI_1\wedge \omega^1-dI_3\wedge \omega^3-I_3K\omega^1\wedge \omega^2-I_1\omega^2\wedge
\omega^3+II_1\omega^1\wedge\omega^3\\
d\theta^2 & =-dK_1\wedge\omega^1-dK_2\wedge\omega^2+IK^2\omega^1\wedge\omega^2+(IK_2-
K_1)\omega^2\wedge\omega^3\\
 & +(2IK_1+I_1K+K_2)\omega^1\wedge\omega^3+K\theta^1\wedge\omega^3+I\theta^2\wedge\omega^3.
\end{split}
\end{equation*}
	\quad Let us remark that the above formulas can be rewritten as
\begin{equation*}
\begin{split}
d\theta^1 & \equiv (-dI_1+I_3K\omega^2-II_1\omega^3)\wedge\omega^1+
(-dI_3-I_1\omega^2)\wedge\omega^3 \quad \mod\ \{I\} \\
d\theta^2 & \equiv \Bigl[-dK_1-IK^2\omega^2-(2IK_1+I_1K+K_2)\omega^3\Bigr]
\wedge\omega^1\\
& +\Bigl[-dK_2-(IK_2-K_1)\omega^3\Bigr]\wedge\omega^2 \quad \mod\ \{I\}.
\end{split}
\end{equation*}
	\quad It follows that we can write
\begin{equation*}
d\theta^a  \equiv \pi^a_i\wedge \omega^i\quad \mod\ \{I\}, 
\end{equation*}
where {\it a}=1,2, {\it i}=1,2,3. \\
	\quad The 1-forms matrix $(\pi^a_i)$ has the following non-vanishing  entries:
\begin{equation}\label{pi_matrix}
\begin{split}
\pi_1^1 & = -dI_1+I_3K\omega^2-II_1\omega^3,\\
\pi_3^1 & = -dI_3-I_1\omega^2,\\
\pi_1^2 & = -dK_1-IK^2\omega^2-(2IK_1+I_1K+K_2)\omega^3\\
\pi_2^2 & = -dK_2-(IK_2-K_1)\omega^3.
\end{split}
\end{equation}
	\quad By putting now 
\begin{equation}\label{pi_vector}
\begin{split}
& \pi^1:=\pi_1^1,\qquad \pi^2:=\pi_3^1,\\
& \pi^3:=\pi_1^2,\qquad \pi^4:=\pi_2^2,
\end{split}
\end{equation}
we obtain that $(I,J)$ is a linear Pfaffian system that lives on the 9 dimensional manifold $\widetilde{\Sigma}$ which has 
the coframing
 \begin{equation}
\{\theta^1,\theta^2,\omega^1,\omega^2,\omega^3,\pi^1,\pi^2,\pi^3,\pi^4
\}
\end{equation}
that is adapted to the filtration 
\begin{equation*}
I\subset J\subset T^*\widetilde{\Sigma}.
\end{equation*}
	\quad Since the apparent torsion was absorbed, we can write
\begin{equation*}
d\theta^a \equiv  A_{\epsilon i}^a\pi^\epsilon \wedge \omega^i \qquad
 \mod\ \{I\},
\end{equation*}
where the non-vanishing entries of  $ A^a_{i\epsilon}$ are
\begin{equation}\label{Lands_tableau}
A_{11}^1=A^1_{23}=A^2_{31}=A^2_{42}=1.
\end{equation}
	\quad The 1-forms $\pi_a^i$ are sections of $T^*\widetilde{\Sigma} / J$, or,
equivalently, they are components of a section of $I^*\otimes J/I$. \\
	\quad From now on, by abuse of notation we will write the structure equations
of the EDS as
\begin{equation*}
\begin{split}
& \theta^a=0\\
& d\theta^a \equiv  A_{\epsilon i}^a\pi^\epsilon \wedge \omega^i \qquad \mod\ \{I\}\\
& \Omega=\omega^1\wedge \omega^2\wedge \omega^3\neq 0.
\end{split}
\end{equation*}
	\quad From (\ref{Lands_tableau}) it follows that 
the tableau $A$ of the linear Pfaffian system $(I,J)$ is given by
\begin{equation*}
A=
\begin{pmatrix}
a & 0 & d \\
b & c & 0
\end{pmatrix},
\end{equation*}
where $a,b,c,d$ are nonzero constants. Therefore, the reduced
characters of the tableau $A$ are $s_1=2$, $s_2=2$, $s_3=0$, and $s_0=$rank$\ I=2$.\\
\quad The symbol $B$ of the linear Pfaffian system $(I,J)$ is then 
\begin{equation*}
B=
\begin{pmatrix}
0 & e & 0 \\
0 & 0 & f
\end{pmatrix},
\end{equation*}
where $e,f$ are nonzero constants. 
\subsection {The integrability conditions}
	\quad Let us denote by $(G_3(T\widetilde{\Sigma}),\pi,\widetilde{\Sigma})$ the Grassmannian of three planes through the 
origin of $T\widetilde{\Sigma}$. Then the dimension of the base manifold and the fiber over a point $p\in 
\widetilde{\Sigma}$ are given by
\begin{equation*}
\dim G_3(T\widetilde{\Sigma}) =27,\qquad \dim G_3(T_p\widetilde{\Sigma}) =18,
\end{equation*}
respectively.\\
	\quad If we denote by $p_i^a$, ($a=1,...,6$, $i=1,2,3$) the local coordinates of the fiber $G_3(T_p
\widetilde{\Sigma})$, then for a $3$-plane $E\in G_3(T_p\widetilde{\Sigma})$, that satisfies the independence condition 
$\omega^1\wedge\omega^2\wedge\omega^3_{\ |E}\neq 0$, by an eventual
relabeling of the coordinates, equations of integral elements of $(I,J)$
are 
\begin{equation}\label{integ_elem}
\begin{split}
& \theta^b=0,\ (b=1,2)\\
& \pi^\epsilon-p^{\epsilon}_i\omega^i=0,
\end{split}
\end{equation}
where $p_i^\epsilon$, ($\epsilon=1,...,4$, $i=1,2,3$) are functions on $G_3(T_p\widetilde{\Sigma})$.\\
	\quad The relations (\ref{integ_elem}) regarded as system of linear equations in
$p_i^\epsilon$ are the {\it first order integrability conditions} of the linear
Pfaffian system $(I,J)$. One can remark that in the most general case,
these equations are over-determined, in the sense that there are more
equations than unknowns. Therefore, in general there is likely for such linear systems to be incompatible.\\ 
	\quad In our case, using the fact that integral elements of $\theta^a=0$ must satisfy 
$d\theta^a=0$ also, then using (\ref{Lands_tableau}) we obtain the solutions
of (\ref{integ_elem}) as follows:
\begin{equation}\label{integ_manif}
\begin{split}
& p_2^1=0,\qquad \qquad \qquad p_3^1=p_1^2,\\
& p_2^2=0,\\
& p_3^3=0,\qquad \qquad \qquad p_2^3=p_1^4,\\
& p_3^4=0,
\end{split}
\end{equation}
the rest of the functions, namely $p^1_1,p^2_1,p^2_3,p^3_1,p^4_1,p^4_2$,
being arbitrary.\\
	\quad One can see that the maximum rank of this system of functions is $d=6$,
and that it is of local rank constant. In other words,  $\mathcal V_3(\mathcal{I},\Omega)$ is a smooth codimension 6 
submanifold of $G_3(T\widetilde{\Sigma})$, where we denoted by $\mathcal V_3(\mathcal{I},\Omega)\subset G_3(T
\widetilde{\Sigma})$ the subbundle of 3-dimensional integral elements of $\mathcal{I}$. \\
\quad Remark that  
$(\mathcal 
V_3(\mathcal{I},\Omega),\widetilde{I})$ is the prolongation of $(\widetilde{\Sigma},\mathcal{I})$, where 
$\mathcal{I}$ is the exterior 
differential system generated by the Pfaffian system $I$. Here, $\widetilde{I}$ is the exterior differential system on 
$\mathcal 
V_3(\mathcal{I},\Omega)$ generated by the Pfaffian system
\begin{equation*}
\widetilde{I}=\{\theta^1,\theta^2,\pi^1-p_1^1\omega^1-p_1^2\omega^3,
\pi^2-p_1^2\omega^1-p_3^2\omega^3,\pi^3-p_1^3\omega^1-p_1^4\omega^2,
\pi^4-p_1^4\omega^1-p_2^4\omega^2\},
\end{equation*} 
i.e. $\widetilde{I}$ is the pullback to $\mathcal 
V_3(\mathcal{I},\Omega)$, by the inclusion $\iota:\mathcal 
V_3(\mathcal{I},\Omega)\to G_3(T_p\widetilde{\Sigma})$,
 of the canonical system on $G_3(T_p\widetilde{\Sigma})$.\\
	\quad Moreover, since the dimension  of the solution space of equations (\ref{integ_manif}) is 6, the Cartan 
involutivity test is satisfied:
\begin{equation*}
s_1+2s_2+3s_3=2+2\cdot2+0=6=d.
\end{equation*}
	\quad Using the Cartan-K\"{a}hler theorem for linear Pfaffian systems
(see \cite{IL2003}, p. 176, \cite{Br et al 1991} for a more general exposition, and the Appendix), we can summarize the 
findings in this section in the following theorem.\\
	\quad{\bf Theorem 4.1.}\\
\quad {\it Assume there exist three 1-forms $(\omega^1, \omega^2,\omega^3)$ on a 9-dimensional manifold
$\widetilde{\Sigma}$ which satisfy the structure equations \eqref{Lands_struct_eq}, where I, K are considered as free coordinates 
on $\widetilde{\Sigma}$, and dI, dK are independent from  $\omega^1, \omega^2,\omega^3$.\\
\quad Then the pair $(I,J)$ is an involutive linear Pfaffian system with independence
condition on $\widetilde{\Sigma}$. Therefore, solving a series of Cauchy problems
yields analytic integral manifolds of $(I,J)$ passing through
$\tilde{u}\in \widetilde{\Sigma}$ that, roughly speaking, 
depend on two functions of two variables.}

	\quad We emphasize the fact that the existence of analytical integral manifolds of $(I,J)$ is guaranteed  only in a 
neighborhood $U\subset \widetilde{\Sigma}$ of $\tilde{u}$. \\
	\quad Therefore, for any point $\tilde{u}\in\widetilde{\Sigma}$ chosen such that $I_1\neq 0$, the existence of an integral 
submanifold of $(I,J)$ passing through this point is guaranteed by Theorem 4.1. This is a non-trivial generalized 
Landsberg surface structure on which the independence condition $\omega^1\wedge\omega^2\wedge\omega^3\neq  0$ 
is satisfied. In other words, this integral submanifold can be realized as the graph of the analytical mapping 
\begin{equation*}
\begin{split}
&\Sigma\to\widetilde{\Sigma},\\
&u \mapsto (u,I(u),K(u),I_1(u),I_3(u),K_1(u),K_2(u))\in \widetilde{\Sigma}.
\end{split}
\end{equation*}
	\quad This proves R. Bryant's prophecy. Unfortunately, these generalized structures are not always amenable, in 
other words, they are not always realizable as Finsler structures on surfaces as will be seen.\\
	\quad {\bf Remark.}\\
\quad If we write the structure equations as
\begin{equation*}
\begin{pmatrix}
d\theta^1\\
d\theta^2
\end{pmatrix}
=
\begin{pmatrix}
\pi^1 & 0  & \pi^2 \\
\pi^3 & \pi^4 & 0 
\end{pmatrix}
\wedge
\begin{pmatrix}
\omega^1\\
\omega^2\\
\omega^3
\end{pmatrix}
,
\end{equation*}
then we can put them in a normal form which reflects
the Cartan test for involutivity. \\
\quad Indeed, if one changes the basis $\{\omega^1,\omega^2,\omega^3\}$ to 
$\{\widetilde{\omega}^1:=\omega^1,\widetilde{\omega}^2:=\omega^2,\widetilde{\omega}^3:=\omega^3-\omega^2\}$, then 
it follows
\begin{equation*}
\begin{split}
& d\theta^1\equiv {\pi}^1\wedge\widetilde{\omega}^1+{\pi}^2\wedge\widetilde{\omega}^2+{\pi}^2\wedge
\widetilde{\omega}^3\\
& d\theta^2\equiv
{\pi}^3\wedge\widetilde{\omega}^1+{\pi}^4\wedge\widetilde{\omega}^2,\qquad\qquad\qquad \mod {I}.
\end{split}
\end{equation*}
	\quad Therefore, in this frame, 
the tableau $A$ of $(I,J)$ is now given by
\begin{equation}
A=
\begin{pmatrix}
a & d & d \\
b & c & 0
\end{pmatrix}.
\end{equation}
	\quad One can now directly verify by visual inspection that, indeed, there are
$s_1=2$ independent 1-forms in the first column of the tableau matrix of
$(I,J)$, $s_1+s_2=4$ independent 1-forms in the first two columns,
and 
$s_1+s_2+s_3=4$ independent 1-forms in the first three columns, i.e. in
the entire matrix. This agrees with Cartan's test for involutivity.  

\subsection{The existence of generalized Landsberg structures on surfaces}
\quad In the present section we are going to generalize our setting and show the existence of the coframes $\omega$ satisfying (\ref{Lands_struct_eq}) together with the scalar functions $I$ and $K$ satisfying (\ref{Lands_Bianchi}), without using any of the assumptions in \S4.1, \S4.2. \\
\quad Let $\Sigma$ be again a 3-manifold, and let $\pi:\mathcal{F}(\Sigma)\to\Sigma$ be its frame bundle, namely
\begin{equation*}
\mathcal{F}(\Sigma)=\{(u,f_u) | f_u:T_u\Sigma\to V  \ \textrm{linear isomorphism}\},
\end{equation*} 
where $V$ is a 3-dimensional real vector space.\\
\quad Let $\eta$ be the tautological $V$-valued 1-form on $\mathcal{F}(\Sigma)$, defined as usual by
\begin{equation}\label{taut_def}
\eta_f(w)=f_u(\pi_*w),
\end{equation}
where $f=(u,f_u)\in \mathcal{F}(\Sigma)$, and $w\in T_f\mathcal{F}(\Sigma)$.\\
\quad It is known that a coframe on the manifold $\mathcal{F}(\Sigma)$ is given by $(\eta^i,\alpha^i_j)$, $i,j=1,2,3$, where $\eta^i$ are the components of the $V$-valued tautological form $\eta$, and $\alpha^i_j$ are the 1-forms on $\mathcal{F}(\Sigma)$ that satisfy the structure equations
\begin{equation}\label{str.eq_1}
d\eta^i=-\alpha^i_j\wedge\eta^j,\qquad i,j=1,2,3.
\end{equation}
\quad Such 1-forms always exist, but without supplementary conditions, they are not unique. These forms are the connection forms of the frame bundle.\\
\quad Here, we choose a "flat type connection form", i.e. 1-forms $\alpha^i_j$ satisfying
\begin{equation}\label{str.eq_2}
d\alpha^i_j=\alpha^i_k\wedge\alpha^k_j,\qquad i,j,k=1,2,3.
\end{equation}
\quad We define next the following (local) trivialization of the frame bundle
\begin{equation}\label{id}
\begin{split}
t:& \ \mathcal{F}(\Sigma)\to \Sigma\times GL(3,\mathbb{R})\\
&f=(u,f_u)\mapsto (u,(f_j^i)),
\end{split}
\end{equation}
where for a coordinate system $(x^1,x^2,x^3)$ on $\Sigma$, and a basis $\{e_1,e_2,e_3\}$ of $V$, $(f_j^i)$ is the representation matrix of the mapping $f_u:T_u\Sigma\to\Sigma$ with respect to the bases $\{\frac{\partial}{\partial x^i}\}$ and $\{e_i\}$.\\
\quad A system of coordinates on $\Sigma\times GL(3,\mathbb{R})$ is given by $(x^i,f^i_j)$, $i,j=1,2,3$, and a coframe on the manifold $\Sigma\times GL(3,\mathbb{R})$ will be $(\omega^i,df^i_j)$, where we put
\begin{equation}\label{omegas}
\omega^i=f_j^idx^j.
\end{equation}
\quad We remark that the tautological 1-forms $\eta=(\eta^i)$, $i=1,2,3$, on $\mathcal{F}(\Sigma)$ correspond to the 1-forms $(\omega^i)$ under the identification (\ref{id}). This can be verified by direct computation checking that the 1-forms $\omega^i$ in (\ref{omegas}) satisfy (\ref{taut_def}).\\
\quad Moreover, if we put
\begin{equation*}
\beta_j^i=d(f_k^i)(f^{-1})_j^k,\qquad i,j,k=1,2,3,
\end{equation*}
then the 1-forms $(\beta_j^i)$ on $\mathcal{F}(\Sigma)$ correspond to the "connection forms" $(\alpha_j^i)$. Indeed, a straightforward computation shows that the $\beta_j^i$'s defined above verify the structure equations (\ref{str.eq_1}), (\ref{str.eq_2}).\\
\quad With these preparations in hand, we move on to the study of the existence of a coframe $\omega$ and the scalars $I$, $K$ on the 3-manifold $\Sigma$ that satisfy (\ref{Lands_struct_eq}) and (\ref{Lands_Bianchi}), respectively. \\
\quad In order to do this, we consider the 18-dimensional manifold 
\begin{equation*}
\widetilde{\Sigma}=\mathcal{F}(\Sigma)\times \mathbb{R}^6
\end{equation*}
with the coframe 
\begin{equation*}
\{\eta^1,\eta^2,\eta^3,(\alpha_j^i)_{i,j=1,2,3},\theta^1,\theta^2,\pi^1,\pi^2,\pi^3,\pi^4\},
\end{equation*}
where $\pi^1,\pi^2,\pi^3,\pi^4$ are the 1-forms in \eqref{pi_matrix}, \eqref{pi_vector}.\\
\quad We consider the 1-forms
\begin{equation}\label{forms1_eta}
\begin{split}
\Theta^1 & =d\eta^1+I\eta^1\wedge\eta^3-\eta^2\wedge\eta^3\\
\Theta^2 & =d\eta^2+\eta^1\wedge\eta^3\\
\Theta^3 & =d\eta^3-K\eta^1\wedge\eta^2
\end{split}
\end{equation}
and
\begin{equation}\label{forms2_eta}
\begin{split}
\theta^1 & =dI-I_1\eta^1-I_3\eta^3\\
\theta^2 & =dK-K_1\eta^1-K_2\eta^2+KI\eta^3,
\end{split}
\end{equation}
obtaining in this way the exterior differential system 
\begin{equation*}
\widetilde{\mathcal I}=\{\Theta^1,\Theta^2,\Theta^3,\theta^1,\theta^2\}
\end{equation*}
with independence condition
\begin{equation*}
\Omega=\eta^1\wedge\eta^2\wedge\eta^3\neq 0.
\end{equation*}

\quad Let us remark that any element $E\in G_3(T\widetilde{\Sigma})$ such that $\Omega |_{E}\neq 0$ 
is defined by 
\begin{equation*}
\begin{split}
\alpha^i_{j\ |E} & =A^i_{jk}(E)\eta^k_{\ |E}\\
\theta^i_{\ |E} & =B_k^i(E)\eta^k_{\ |E}\\
\pi^i_{\ |E} & =C_k^i(E)\eta^k_{\ |E},
\end{split}
\end{equation*}
where $(A^i_{jk})_{i,j,k=1,2,3}$, $(B_k^i)_{i=1,2;k=1,2,3}$, 
$(C_k^i)_{i=1,2,3,4;k=1,2,3}$ 
are smooth functions on $G_3(T\widetilde{\Sigma},\Omega)$. 
\quad In other words, $(A^i_{jk},B_k^i,C_k^i)$ are the fiber coordinates of the fibration
$G_3(T\widetilde{\Sigma},\Omega)\to \widetilde{\Sigma}$. This fiber is 45-dimensional.\\
\quad However, due to the identification (\ref{id}) and the discussion above, we can consider the local coordinates
\begin{equation*}
(x^1,x^2,x^3,(f_j^i)_{i,j=1,2,3},I,K,I_1,I_3,K_1,K_2)\in\widetilde{\Sigma}
\end{equation*}
on the 18-dimensional manifold $\mathcal{F}(\Sigma)\times\mathbb{R}^6$ and identify the 1-forms $\eta^i$ with 
$\omega^i$ given in (\ref{omegas}). Since the settings are equivalent, for simplicity, we will work in these coordinates instead of the general case described at the beginning of this subsection. \\
\quad It follows that the 1-forms (\ref{forms1_eta}), (\ref{forms2_eta}) of the exterior differential system 
$\widetilde{\mathcal{I}}$, can be written as
\begin{equation}\label{forms1_omega}
\begin{split}
\Theta^1 & =d\omega^1+I\omega^1\wedge\omega^3-\omega^2\wedge\omega^3\\
\Theta^2 & =d\omega^2+\omega^1\wedge\omega^3\\
\Theta^3 & =d\omega^3-K\omega^1\wedge\omega^2
\end{split}
\end{equation}
and
\begin{equation}\label{forms2_omega}
\begin{split}
\theta^1 & =dI-I_1\omega^1-I_3\omega^3\\
\theta^2 & =dK-K_1\omega^2-K_2\omega^2+KI\omega^3,
\end{split}
\end{equation}
with independence condition
\begin{equation*}
\Omega=\omega^1\wedge\omega^2\wedge\omega^3\neq 0,
\end{equation*}
where $\omega$'s are given by (\ref{omegas}).\\
\quad The integral manifolds of $(\widetilde{\mathcal{I}},\Omega)$ will consist of the coframe 
$\{\omega^1,\omega^2,\omega^3\}$, and the functions $(I,K,I_1,I_3,K_1,K_2)$ on $\Sigma$. The projection of such integral manifold to $\Sigma$ gives a generalized Landsberg structure $(\Sigma,\omega)$.\\
\quad Let us remark that the situation is now quite different from the one in \S4.1. The $\Theta$'s are 2-forms, while $\theta$'s are 1-forms, so the exterior differential system  $(\widetilde{\mathcal{I}},\Omega)$ is not a linear Pfaffian system, and therefore we cannot apply the Cartan-K\"ahler theorem for linear Pfaffian systems as we did previously. Even there are more general versions of the  Cartan-K\"ahler theorem, the strategy we adopt here is to prolong   
$\widetilde{\mathcal{I}}$ in order to obtain a linear Pfaffian system (for details see \cite{IL2003}, p. 177).\\
\quad Let us consider the prolongation $\mathcal{V}(\widetilde{\mathcal{I}},\Omega)\subset G_3(T\widetilde{\Sigma})$
over $\widetilde{\Sigma}$, with the fiber inhomogeneous Grassmannian coordinates
 $\Bigl((p^i_j)_{i=1,2;j=1,2,3}, (p_{jk}^i)_{i,j,k=1,2,3}$, \\
 $(q_k^i)_{i=1,2,3,4;k=1,2,3}\Bigr)$, such that
\begin{equation*}
\begin{split}
 \theta^i_{|_E} &  =p_k^i(E)dx^k_{|_E}\\
 {df_j^i}_{|_E}  & =p_{jk}^i(E)dx^k_{|_E}\\
 \pi^i_{|_E}  & =q_k^i(E)dx^k_{|_E},
\end{split}
\end{equation*} 
for any integral element $E$.\\
\quad Then, the equations
\begin{equation*}
\begin{split}
& \theta^i=d\theta^i=0,\qquad i=1,2,\\
& \Theta^j=d\Theta^j=0,\qquad j=1,2,3,
\end{split}
\end{equation*}
will give the defining equations of the prolongation $\mathcal{V}(\widetilde{\mathcal{I}},\Omega)$.\\
\quad As concrete computation, we remark first that $\theta^i=0$ will imply $p_j^i=0$, so these functions will not appear in out analysis. A similar computation as in \S4.1 shows that the structure equations for $\theta$'s are
\begin{equation*}
\begin{split}
& d\theta^1\equiv \pi^1\wedge \omega^1+\pi^2\wedge \omega^2\quad \mod \{\theta,\Theta\} \\
& d\theta^2\equiv \pi^3\wedge \omega^1+\pi^4\wedge \omega^3.
\end{split}
\end{equation*}
These equations will give some of the $q_j^i$'s.\\
\quad The equations $\Theta^i\equiv 0$ $\mod \{\theta,\Theta\}$ will give some of the $p_{jk}^i$. The rest of the equations $d\Theta^i\equiv 0$ $\mod \{\theta,\Theta\}$ will be satisfied due to some Bianchi identities, so they will give no further conditions.\\
\quad In this way, we obtain the linear Pfaffian $\widetilde{\widetilde{\mathcal{I}}}$ 
on $\mathcal{V}(\widetilde{I},\Omega)$
generated by the 1-forms
\begin{equation}\label{large_lin_pfaff}
\{\theta^1,\theta^2,({\Theta^i_j})_{i,j=1,2,3},\Pi^1,\Pi^2,\Pi^3,\Pi^4\},
\end{equation}
where 
\begin{equation*}
\begin{split}
& \Theta^i_j=df^i_j-p_{jk}^idx^k,\quad i,j,k=1,2,3,\\
& \Pi^i=\pi^i-q^i_kdx^k,\quad i=1,2,3,4, \quad k=1,2,3
\end{split}
\end{equation*}
and we will study its involutivity by means of Cartan-K\"ahler theory as we did in \S4.1, \S4.2.\\
\quad It is easy to see that putting the conditions $d\theta^i=0$, $i=1,2$ it results 6 relations with 12 unknown functions $(q_j^i)_{i=1,2,3,4;\ j=1,2,3}$. We solve $q_3^1$, $q_2^2$, $q_3^2$ in terms of $q_1^1$, $q_2^1$, $q_1^2$, and  
$q_3^3$, $q_2^4$, $q_3^4$ in terms of $q_1^3$, $q_2^3$, $q_1^4$. It follows
\begin{equation*}
\begin{split}
& \quad q_2^2=\frac{1}{f_1^3}\Bigl(q_1^1f_2^1-q_2^1f_1^1+q^2_1f_1^3\Bigr),\\
& \begin{pmatrix}
q_3^1\\
q_3^2
\end{pmatrix}=
\begin{pmatrix}
-f_2^1 & -f_2^3 \\
f_1^1 & f_1^3
\end{pmatrix}^{-1}
\begin{pmatrix}
-q_2^1f_3^1 & -q_2^2f_3^2 \\
q_1^1f_3^1 & q_1^2f_3^3
\end{pmatrix},
\end{split}
\end{equation*}
and
\begin{equation*}
\begin{split}
& \quad q_2^4=\frac{1}{f_1^2}\Bigl(q_3^1f_2^1-q_2^3f_1^1+q^4_1f_1^2\Bigr),\\
& \begin{pmatrix}
q_3^3\\
q_3^4
\end{pmatrix}=
\begin{pmatrix}
-f_2^1 & -f_2^2 \\
f_1^1 & f_1^2
\end{pmatrix}^{-1}
\begin{pmatrix}
-q_2^1f_3^1 & -q_2^4f_3^2 \\
q_1^3f_3^1 & q_1^4f_3^2
\end{pmatrix}.
\end{split}
\end{equation*}
\quad In the same way, from $\Theta^i_j=0$, $i,j=1,2,3$, we obtain 9 relations with 27 unknown functions $(p_{jk}^i)$, $i,j,k=1,2,3$. Solving 9 of them, we obtain
\begin{equation*}
\begin{split}
& p_{31}^1=p_{13}^1-I(f_3^1f_1^3-f_1^1f_3^3)+(f_3^2f_1^3-f_1^2f_3^3)\\
& p_{12}^1=p_{21}^1-I(f_2^1f_1^3-f_1^1f_2^3)+(f_2^2f_1^3-f_1^2f_2^3)\\
& p_{23}^1=p_{32}^1-I(f_3^1f_2^3-f_2^1f_3^3)+(f_3^2f_2^3-f_2^2f_3^3),
\end{split}
\end{equation*}
\begin{equation*}
\begin{split}
& p_{12}^2=p_{21}^2-(f_1^1f_2^3-f_2^1f_1^3)\\
& p_{23}^2=p_{32}^2-(f_2^1f_3^3-f_3^1f_2^3)\\
& p_{31}^3=p_{13}^3-(f_3^1f_1^3-f_1^1f_3^3),
\end{split}
\end{equation*}
\begin{equation*}
\begin{split}
& p_{12}^3=p_{21}^3+K(f_1^1f_2^2-f_2^1f_1^2)\\
& p_{23}^3=p_{32}^3+K(f_2^1f_3^2-f_3^1f_2^2)\\
& p_{31}^3=p_{13}^3+K(f_3^1f_1^2-f_1^1f_3^2).
\end{split}
\end{equation*}
\quad Using these relations we study the involutivity of the linear Pfaffian (\ref{large_lin_pfaff}). By similar computations as in \S4.1, \S4.2 we obtain that the structure equations of (\ref{large_lin_pfaff}) are given by
\begin{equation}\label{large_str.eq}
d\begin{pmatrix} \theta^1 \\ \theta^2 \\ \Theta^1_1 \\ \Theta^1_2 \\ \Theta^1_3 \\  \Theta^2_1 \\ \Theta^2_2 \\ \Theta^2_3 \\  \Theta^3_1 \\ \Theta^3_2 \\ \Theta^3_3 \\ \Pi^1 \\ \Pi^2 \\ \Pi^3 \\ \Pi^4
\end{pmatrix}\equiv
\begin{pmatrix}
0 & 0 & 0 \\
0 & 0 & 0 \\
\rho^1 & \rho^2 & \rho^3\\
\rho^2 & \rho^4 & \rho^5\\
\rho^3 & \rho^5 & \rho^6\\
\rho^7 & \rho^8 & \rho^9\\
\rho^8 & \rho^{10} & \rho^{11}\\
\rho^9 & \rho^{11} & \rho^{12}\\
\rho^{13} & \rho^{14} & \rho^{15}\\
\rho^{14} & \rho^{16} & \rho^{17}\\
\rho^{15} & \rho^{17} & \rho^{18}\\
\rho^{19} & \rho^{20} & \Phi^1\\
\rho^{21} & \Phi^2 & \Phi^3\\
\rho^{22} & \rho^{23} & \Phi^4\\
\rho^{24} & \Phi^5 & \Phi^6
\end{pmatrix}
\begin{pmatrix}
dx^1 \\ dx^2 \\ dx^3
\end{pmatrix}\qquad \mod\{\widetilde{\widetilde{\mathcal{I}}}\},
\end{equation}
where $\rho^i$, $i=1,\dots,24$ are 1-forms on $\mathcal{V}(\widetilde{\mathcal{I}},\Omega)$, 
linearly independent from the 1-forms in (\ref{large_lin_pfaff}), and $\Phi^j$, $j=1,\dots,6$ are linear combinations of the $\rho$'s. \\
\quad It means that the apparent torsion can be absorbed. It also can be checked that the space of integral elements at each point has dimension 38.\\
\quad On the other hand, the reduced characters of the tableau corresponding to (\ref{large_str.eq}) are
\begin{equation*}
s_1=13,\quad s_2=8,\quad s_3=3,
\end{equation*}
and Cartan test's for involutivity reads
\begin{equation*}
s_1+2s_2+3s_3=38.
\end{equation*}
Therefore the Pfaffian system (\ref{large_lin_pfaff}) is involutive.\\
\quad Putting all these together, and assuming that $\Sigma$ and  $\alpha,\eta$ are 
analytic, from Cartan-K\"ahler theory we obtain

\quad {\bf Theorem 4.2.}\\
\quad {\it The linear Pfaffian prolongation 
$(\mathcal{V}(\widetilde{\mathcal{I}},\Omega),\widetilde{\widetilde{\mathcal{I}}})$ of the exterior differential system $\widetilde{\mathcal{I}}$ on $\widetilde{\Sigma}$ is involutive. Moreover, the analytical 
integral manifolds of $\widetilde{\widetilde{\mathcal{I}}}$ depend on 3 functions of 3 variables.
}\\
\quad Since the projection of an integral manifold of the prolongation $\widetilde{\widetilde{\mathcal{I}}}$ to $\widetilde{\Sigma}$ is also an integral manifold of $\widetilde{\mathcal{I}}$, it follows

\quad {\bf Corollary 4.3.}\\
\quad {\it There exist non-trivial generalized Landsberg structures on a 3-manifold $\Sigma$.}

\quad The non-triviality of the integral manifolds can be obtained by choosing an appropriate initial value. See the discussion at the end of \S4.2.\\
\quad {\bf Remark.}\\
\quad We point out that the degree of freedom of the integral manifolds of $\widetilde{\widetilde{\mathcal{I}}}$ does not equal the degree of freedom of the scalar functions $I$ and $K$. The reason is that the 3 functions of 3 variables obtained in Theorem 4.2 include the degree of freedom of the coframe $(\omega^1,\omega^2,\omega^3)$ as well.


\section{The local amenability of generalized Landsberg structures on surfaces}

\quad The notion of amenability given in Definition 2.2 has the following local version 

\quad  {\bf Definition 5.1.} The generalized Finsler structure $(\Sigma,\omega)$ is 
called {\it locally} {\it amenable} if for any point $u\in \Sigma$, there exists an open neighborhood $U\subset \Sigma$ of $u$ to which  $(\Sigma,\omega)$ restricts to be amenable, i.e. $(U,\omega_{|U})$ is amenable in the sense of Definition 2.2.

\quad We can now formulate a local version of the Theorem 2.1.

 {\bf Theorem 5.1.}\quad {\it Let $(\Sigma,\omega)$ be a generalized Finsler 
structure. Then the following two conditions always hold good.

(2)' $(\Sigma,\omega)$ is locally amenable.

(3)' The mapping $\nu:U\to T(\widetilde{U})$ is a smooth embedding, where 
$\widetilde{U}$ is the leaf space of the foliation 
$\{\omega^1_{|U}=0,\omega^2_{|U}=0\}$.
}

{\it Proof.} The proof is quite straightforward. Remark first that the differential system $\{\omega^1=0,\omega^2=0\}$ is completely integrable. Indeed, the structure equations (\ref{finsler_struct_eq}) of a generalized Finsler structure show that 
\begin{equation*}
\begin{split}
& d\omega^1 \equiv 0\qquad \qquad \mod\{\omega^1,\omega^2\}.\\
& d\omega^2 \equiv 0
\end{split}
\end{equation*}
\quad It follows from Frobenius theorem that for any point $u\in \Sigma$, there exists an open neighborhood $U\subset\Sigma$ of $u$ such that the leaf space of the foliation $\{\omega^1_{|U}=0,\omega^2_{|U}=0\}$ is a differentiable manifold, say $\widetilde{U}$, such that the canonical projection $\pi:U\to \widetilde{U}$ is a smooth submersion.

\quad From here we see immediately that $\nu:U\to T(\widetilde{U})$ is a smooth embedding.
\begin{flushright}
Q. E. D.
\end{flushright}

\quad We point out that the condition (1) in Theorem 2.1 is not necessarily true for this $U$.

\quad Indeed, imagine for a moment the case when  the generalized Finsler structure $(\Sigma,\omega)$ satisfies all the conditions in Theorem 2.1, i.e. it is a classical Finsler structure on a differentiable surface $M$ such that $\pi:\Sigma\to M$ is a smooth submersion. In this case, even though if we restrict ourselves to a small neighborhood 
$\widetilde{U}\subset M$, the fibers $\Sigma_x$ over $x\in \widetilde{U}$ are not changed in any way, they remain diffeomorphic to $S^1$ when we shrink the base manifold $M$.

\quad This situation changes dramatically when we are working with a local generalized structure on $\Sigma$. Considering the neighborhood $U \subset \Sigma$ as given by the Frobenius theorem, the fibers are also {\it cut off}. The situation is similar with taking a neighborhood of a point on the surface of the sphere 
$S^2$, for example. In general, the great circles will have only some open arcs contained in this neighborhood, and there is no reason for these arcs to be compact.

\quad Hence, the local conditions in Theorem 5.1 are not enough for $(\Sigma, \omega)$ to be classical Finsler structure on $\widetilde{U}$. 

\quad Therefore, we have

{\bf Corollary 5.2.}\quad {\it Let $(\Sigma,\omega)$ be a generalized Finsler structure and let $U\subset \Sigma$ be the neighborhood given in Theorem 5.1, where (2)' and (3)' are satisfied.
\quad Then $(U,\omega_{|U})$ satisfies (1) in Theorem 2.1 if and only if it is a classical Finsler structure on  $\widetilde{U}$. 
}

\quad In conclusion, recall that we have proved the existence of non-trivial generalized Landsberg surfaces in Theorem 4.1. In other words, the 
Cartan--K\"ahler theorem assures us that there exists a neighborhood $U\subset \Sigma$ such that $(U,\omega_{|U})$ is a non-trivial generalized Landsberg surface. 

\quad On the other hand, since the differential system $(U,\omega^1,\omega^2)$ is completely integrable, from the discussion above it follows that, on a possible smaller open set on $\Sigma$, there exists the local coordinate system $u=(x,y,p)$ such that the leaf space of the foliation $\{\omega^1=0,\omega^2=0\}$ is a differentiable manifold. 

\quad We can therefore conclude that for a small enough $\varepsilon>0$, there exist amenable non-trivial generalized Landsberg structures $(U,\omega)$, depending on two functions of two variables, over an open disk $D=\{(x,y):x^2+y^2<\varepsilon\}\subset \widetilde{U}$ in the plane.

\quad Finally, we emphasize that these non-trivial Landsberg generalized structures do not necessarily satisfy the condition (1) in Theorem 2.1, so they are not necessarily classical Finsler structures.


\section{A special coframing}


\quad For a nowhere vanishing smooth function $m$ on $\Sigma$, we define the 1-forms
\begin{equation}\label{coframe_change}
 \begin{split}
  \theta^1 & =m\omega^2\\
  \theta^2 & =\omega^3\\
  \theta^3 & =m\omega^1+m_3\omega^2,
 \end{split}
\end{equation}
where the subscripts represent the directional derivatives with respect to the generalized Landsberg
coframe $(\omega^1,\omega^2,\omega^3)$.
\\
\quad Remark that 
\begin{equation*}
 \theta^1\wedge\theta^2\wedge\theta^3=m^2\omega^1\wedge\omega^2\wedge\omega^3,
\end{equation*}
therefore $\{\theta^1,\theta^2,\theta^3\}$ is a coframe on $\Sigma$ provided $m$ is nowhere vanishing smooth function on $
\Sigma$.\\
\quad An easy linear algebra exercise will show that we have
\begin{equation}\label{frame_change}
 \begin{split}
f_1 & =-\frac{m_3}{m}\hat{e}_1+\frac{1}{m}\hat{e}_2\\
f_2 & =\hat{e}_3\\
f_3 & =\frac{1}{m}\hat{e}_1,
 \end{split}
\end{equation}
where we have denoted by $\{f_1,f_2,f_3\}$ and $\{\hat{e}_1,\hat{e}_2,\hat{e}_3\}$ the dual frames of 
$\{\theta^1,\theta^2,\theta^3\}$ and $\{\omega^1,\omega^2,\omega^3\}$, respectively.\\
\quad We would like to impose conditions on the function $m$ and the invariants $I$, $K$ such that the new coframe 
$\theta=\{\theta^1,\theta^2,\theta^3\}$ satisfies the structure equations
\begin{equation}\label{k_struct_eq}
 \begin{split}
  d\theta^1 & =\theta^2\wedge\theta^3\\
  d\theta^2 & =\theta^3\wedge\theta^1\\
  d\theta^3 & = k\theta^1\wedge\theta^2,
 \end{split}
\end{equation}
where $k$ is a smooth function on $\Sigma$ to be determined (one can see that from the third structure equation of the 
coframe $\theta$ that $dk\wedge \theta^1\wedge \theta^2=0$, therefore the directional derivative of $k$ with respect to $
\theta^3$  must vanish). This is a so-called {\it K-Cartan structure} (see \cite{GG2002}). \\
\quad Straightforward computations show that
\begin{equation*}
d\theta^1  =\theta^2\wedge\theta^3
\end{equation*}
holds if and only if
\begin{equation*}
 m_1=0.
\end{equation*}
This is our first condition on $m$.\\
\quad It also follows that 
\begin{equation}
 I=-2\ \frac{m_3}{m},\quad K=m^2.
\end{equation}\\
\quad In this case, we obtain
 \begin{equation}
k=1-\frac{m_{33}}{m}.
\end{equation}
%
%
%
\quad Remark that the {\it Landsberg condition} reads
\begin{equation*}
 I_2=0 \Longleftrightarrow m_{32}=\frac{m_2m_3}{m},
\end{equation*}
and the {non-triviality conditions}
\begin{equation*}
\begin{split}
I_1 & \neq 0 \Longleftrightarrow m_2\neq 0 \\
I_3 & \neq 0 \Longleftrightarrow mm_{33}-(m_3)^2\neq 0\\
K_2 & \neq 0 \Longleftrightarrow I_1  \neq 0\\
K_3 & \neq 0 \Longleftrightarrow m_3 \neq 0. 
\end{split}
\end{equation*}
%
%
\quad We obtain therefore the following

{\bf Proposition 6.1.} \\
{\it Let $(\Sigma, \omega)$ be a generalized Landsberg structure on the 3-manifold $\Sigma$ and let $m:\Sigma \to \mathbb{R}$ be a smooth nowhere vanishing function satisfying the conditions
\begin{enumerate}
\item the direction invariance condition
\begin{equation}
m_1=0
\end{equation}
\item the Landsberg condition
\begin{equation}
m_{23}=\frac{m_2 m_3}{m}.
\end{equation}
\end{enumerate}
Then $\theta=\{\theta^1, \theta^2, \theta^3\}$, with the $\theta^i$'s given in (\ref{coframe_change}), is a coframe on the 3-manifold $\Sigma$ that satisfies the structure equations (\ref{k_struct_eq}) with 

(3) the curvature condition
\begin{equation}\label{curv_cond}
k=1-\frac{m_{33}}{m}.
\end{equation}

 }
\quad Remark that in this case, besides the conditions in the proposition above, the function $m$ will satisfy the Ricci type identities 
\begin{equation*}
\begin{split}
& m_{21}=-m^2m_3\\
& m_{23}-m_{32}=0\\
& m_{31}=m_2.
\end{split}
\end{equation*}
\quad Conversely, we can start with a coframe $\theta=\{\theta^1, \theta^2, \theta^3\}$ on the 3-manifold $\Sigma$ that satisfies the structure equations (\ref{k_struct_eq}) for a function $k:\Sigma\to \mathbb{R}$ such that $k_{\theta 3}=0$. Here, we denote by 
$h_{\theta i}$ the directional derivatives of a smooth function $h$ with respect to the coframe $\theta$, i.e. $dh=h_{\theta 1}\theta^1+h_{\theta 2}\theta^2+h_{\theta 3}\theta^3$.
Making use of a nowhere vanishing smooth function $m:\Sigma\to \mathbb{R}$, we can construct the 1-forms 
\begin{equation}\label{coframe_change_omega}
\begin{split}
&\omega^1=\frac{1}{m}(\theta^3-\frac{m_{\theta 2}}{m}\theta^1)\\
&\omega^2=\frac{1}{m}\theta^1\\
&\omega^3=\theta^2.
\end{split}
\end{equation}

\quad By a simple straightforward computation we obtain

{\bf Proposition 6.2.}

{\it Let $\theta=\{\theta^1, \theta^2, \theta^3\}$ be a coframe on the 3-manifold $\Sigma$ that satisfies the structure equations (\ref{k_struct_eq}) for a smooth function 
$k:\Sigma \to \mathbb{R}$, and let
$m:\Sigma\to \mathbb{R}$ be a nowhere vanishing smooth function that satisfies the conditions
\begin{enumerate}
\item the direction invariance condition
\begin{equation}
m_{\theta 3}=0,
\end{equation}
\item the Landsberg condition
\begin{equation}\label{lands_cond_theta}
(L)\qquad m_{\theta 21}=0,
\end{equation}
\item the curvature condition
\begin{equation}\label{curv_cond_theta}
(C)\qquad \frac{m_{\theta 22}}{m}=1-k.
\end{equation}
Then $\omega=\{\omega^1,\omega^2,\omega^3\}$, with the $\omega^i$'s given in (\ref{coframe_change_omega}), is a generalized Landsberg structure on the 3-manifold $\Sigma$ with the invariants
\begin{equation}
I=-2\frac{m_{\theta 2}}{m},\quad K=m^2.
\end{equation}
\end{enumerate}
}
\quad In this case, the Ricci type equations for $m$ in the coframe $\theta^1$, $\theta^2$, $\theta^3$ are
\begin{equation}\label{Ricci_cond_theta}
\begin{split}
& m_{\theta 12}=m_{\theta 21}=0\\
& m_{\theta 13}=-m_{\theta 2}\\
& m_{\theta 23}=m_{\theta 1}.
\end{split}
\end{equation}

\quad{\bf Remarks.}
\begin{enumerate}
\item Let $(\Sigma,\omega)$ be a generalized Landsberg structure, and suppose that $U\subset \Sigma$ is an open set 
where the foliation  
\begin{equation*}
 \mathcal R = \{\omega^2=0,\omega^3=0\}
\end{equation*}
 is amenable, i.e. the leaf space $\Lambda$ of integral curves of $\hat e_1$ in $U$ is a differentiable manifold, and
\begin{equation*}
l:U\to \Lambda
\end{equation*}
 is a smooth submersion. Then $\theta^1$, $\theta^2$ can be regarded as the tautological 1-forms of the frame bundle and $\theta^3$ as the Levi-Civita connection of the Riemannian manifold 
 $\Lambda$. The function $k$ plays the role of the Gauss curvature.
 \item The indicatrix foliation $\mathcal Q:\{\omega^1=0,\omega^2=0\}$ of the generalized Landsberg structure 
$\{\omega^1,\omega^2,\omega^3\}$ coincides with the geodesic foliation $\mathcal P:\{\theta^1=0,\theta^3=0\}$ of the new 
coframe 
$\{\theta^1,\theta^2,\theta^3\}$ on $\Sigma$.
\item  The normal foliation $\mathcal R:\{\omega^2=0,\omega^3=0\}$ of the generalized Landsberg structure 
$\{\omega^1,\omega^2,\omega^3\}$ coincides with the indicatrix foliation $\mathcal Q:\{\theta^1=0,\theta^2=0\}$ of the 
coframe 
$\{\theta^1,\theta^2,\theta^3\}$ on $\Sigma$.
\item In the case when the generalized Landsberg structure $\{\omega^1,\omega^2,\omega^3\}$ is realizable as a 
classical Finsler structure $(M,F)$ on a certain 2-dimensional differentiable manifold $M$ such that $\pi:\Sigma \to M$ is 
its indicatrix bundle, then the leaves of the normal foliation $\mathcal R:\{\omega^2=0,\omega^3=0\}$ are the (normal) 
lifts of some paths on $M$ called $N$-parallel or $N$-extremal curves. The geometric meaning of such curves $\gamma:
[a,b]\to M$ is that the normal vector field $N(t)$ along $\gamma(t)$, defined by $g_N(N,T)=0$, is parallel 
along $\gamma$. 
Here $T(t)$ is the tangent vector field to the curve $\gamma$, and $g$ is the Riemannian metric induced by the 
Finslerian structure in each tangent plane $T_xM$. It is also known that the $N$-parallels $\gamma$ are solutions of a 
second order differential equation on $M$ and the solution of this SODE is uniquely determined by some initial conditions 
$(x_0,Y_0)\in TM$ (see \cite{ISS2009} for details).  
\end{enumerate}

\section{The geometry of quotient space $\Lambda$}
\subsection{The setting}
\quad In the light of our discussion in \S6, we can conclude that if $U\subset \Sigma$ is an open set where the normal foliation 
$ \mathcal R = \{\omega^2=0,\omega^3=0\}$ 
is amenable, i.e. the leaf space $\Lambda$ of integral curves of $\hat e_1$ in $U$ is a 
differentiable manifold,
 $l:U\to \Lambda$ is a smooth submersion, and $m$ is a smooth function on $\Sigma$ that satisfies the conditions in Proposition 6.1., 
 then there exist
\begin{enumerate}
\item a  quadratic form $g$ on $\Lambda$ such that $l^*(g)=m^2(\omega^2)^2+(\omega^3)^2$;
 \item a 2-form $dA$ on $\Lambda$ such that 
$l^*(dA)=m\omega^2\wedge\omega^3$;
\item a smooth function $\bar{m}$ on $\Lambda$ such that $l^*(\bar{m})=m$.
\end{enumerate}
\quad We can construct now a $g$-orthonormal coframe $\eta^1$, $\eta^2$ on $\Lambda$ (it may be only locally defined), i.e. there exist two 1-forms $
\eta^1$, $\eta^2$ on $\Lambda$, such that 
\begin{equation*}
g=(\eta^1)^2+(\eta^2)^2,\qquad dA=\eta^1\wedge \eta^2>0.
\end{equation*}
\quad This is equivalent with giving a smooth section $s$ of the orthonormal frame bundle $\nu:\mathcal F(\Lambda)
\longrightarrow \Lambda$, i.e. a {\it first order adapted lift} to the geometry of the Riemannian manifold $(\Lambda,g)$.\\
\quad If we denote by $\{e_1,e_2\}$ the dual frame of $\{\eta^1,\eta^2\}$ it follows that $\{e_{1\ |z},e_{2\ |z}\}$ is a 
$g$-orthonormal basis of $T_z\Lambda$, and $(z,e_{1\ |z},e_{2\ |z})\in \mathcal F(\Lambda)$ is a frame on the manifold $
\Lambda$ at each point $z\in \Lambda$. \\
\quad There exist two smooth functions, say $a$ and $b$, on  $\Lambda$ 
such that
\begin{equation*}
\begin{split}
& d\eta^1=a\eta^1\wedge\eta^2\\
& d\eta^2=b\eta^1\wedge \eta^2.
\end{split}
\end{equation*}
\quad By straightforward computation, it also follows that there exists a 1-form, say $\eta^3$, on $\Lambda$, such that
\begin{equation*}
\begin{split}
& d\eta^1=\eta^2\wedge\eta^3\\
& d\eta^2=\eta^3\wedge \eta^1,
\end{split}
\end{equation*}
and therefore we must have
\begin{equation*}
\eta^3=-a\eta^1-b\eta^2.
\end{equation*}
\quad One can easily check that if $\{\widetilde{\eta}^1,\widetilde{\eta}^2\}$ is another $g$-orthonormal frame, then it 
follows $d\widetilde{\eta}^3=d\eta^3$.\\
\quad By straightforward computation we obtain further
\begin{equation*}
d\eta^3=R\eta^1\wedge\eta^2,
\end{equation*}
where $R=a_2-a^2-b_1-b^2$, where $a_i$, $b_i$ means directional derivatives with respect to the coframe 
$\{\eta^1,\eta^2\}$. \\
\quad One can easily see that for another $g$-orthonormal frame  $\{\widetilde{\eta}^1,\widetilde{\eta}^2\}$, the function 
$R$ remains unchanged, and therefore it depends only on $g$. \\
\quad Let us denote
\begin{equation*}
\begin{split}
s(z)=(z,f_z)
\end{split}
\end{equation*}
a local section of $\nu:\mathcal F(\Lambda)\to \Lambda$.\\
\quad It is then known that on $\mathcal F(\Lambda)$ there are tautological 1-forms
\begin{equation*}
\alpha^i_f\in T_f^* \mathcal F(\Lambda),\quad \alpha^i_f:=\eta^i(\nu_*w),
\end{equation*}
where $w\in  T_f \mathcal F(\Lambda)$, and $i\in\{1,2\}$, such that
\begin{equation*}
(\nu^*_f(\eta^1),\nu^*_f(\eta^2))=(\alpha^1_f,\alpha^2_f)
\end{equation*}
gives a basis of semibasic forms on $\mathcal F(\Lambda)$.\\
\quad Consider now the $g$-orthonormal frame bundle 
$\nu:\mathcal F_{\textrm{on}}(\Lambda)\to \Lambda$ with its tautological 1-forms $\{\alpha^1,\alpha^2\}$.\\
\quad If $s:\Lambda\to \mathcal F_{\textrm{on}}(\Lambda)$ is a smooth (local) section, then 
\begin{equation*}
\begin{split}
& \eta^1=s^*(\alpha^1)\\
& \eta^2=s^*(\alpha^2)
\end{split}
\end{equation*}
is a local coframe on $\Lambda$ such that
\begin{equation*}
g=(\eta^1)^2+(\eta^2)^2.
\end{equation*}
\quad Recall that the {\it ``downstairs''} Fundamental Lemma of Riemannian geometry 
tells us that there exists a unique 1-form $\eta^3$ on $\Lambda$ such that 
\begin{equation*}
\begin{split}
& s^*(d\alpha^1)=s^*(\alpha^2)\wedge\eta^3\\
& s^*(d\alpha^2)=\eta^3\wedge s^*(\alpha^1)\\
& d\eta^3=Rs^*(\alpha^1)\wedge s^*(\alpha^2),
\end{split}
\end{equation*}
where $R:\Lambda\to \mathbb R$ is the Gauss curvature of the Riemannian surface $(\Lambda,g)$. These are the so-
called {\it ``downstairs''} structure equations of the Riemannian metric $g$ on $\Lambda$.\\
\quad We also recall the  {\it ``upstairs''} Fundamental Lemma of Riemannian geometry
that states that it must exist a unique 1-form $\alpha^3$ on 
$\mathcal F(\Lambda)$ such that
\begin{equation*}
\begin{split}
& d\alpha^1=\alpha^2\wedge \alpha^3\\
& d\alpha^2=\alpha^3\wedge\alpha^1\\
& d\alpha^3=k\alpha^1\wedge\alpha^2,
\end{split}
\end{equation*}
where $k:\mathcal F(\Lambda)\to \mathbb R$ is the Gauss curvature {\it ``upstairs''}. In our setting it must satisfy the curvature condition \eqref{curv_cond}.
It follows that $R=s^*k$. These are 
the  {\it ``upstairs''} structure equations of the Riemannian metric $g$ on $\Lambda$. One can also see that on $\Lambda
$ we have
\begin{equation*}
(\eta^1,\eta^2,\eta^3)=s^*(\alpha^1,\alpha^2,\alpha^3).
\end{equation*}
%
\quad {\it Example 7.1.}\\
\quad Let us consider a flat Riemannian metric $\widetilde g$ on $\Lambda$, i.e. $\widetilde R=0$. It follows that there 
exist local coordinates $z=(z^1,z^2)$ on $\Lambda$,
 such that
\begin{equation*}
 \widetilde \eta^1=dz^1,\qquad \widetilde \eta^2=dz^2,
\end{equation*}
and therefore
$a=0$, $b=0$ because $d\widetilde \eta^1=0$, $d\widetilde \eta^2=0$.\\
\quad It follows $\widetilde \eta^3=0$ as well as $R=0$.\\
\quad We construct now the coframe $(z; dz^1,dz^2)$ on $\Lambda$ and its oriented orthonormal frame bundle 
$\nu:\widetilde{\mathcal F}_{\textrm{on}}(\Lambda)\to \Lambda$
with respect to the Riemannian metric
\begin{equation*}
\tilde{g}=(dz^1)^2+(dz^2)^2.
\end{equation*}\\
\quad In this case, the tautological 1-forms on $\widetilde{\mathcal F}_{\textrm{on}}(\Lambda)$ will have the normal form
\begin{equation*}
\begin{split}
& \widetilde \alpha^1=\cos(t) dz^1-\sin(t) dz^2\\
& \widetilde \alpha^2=+\sin(t) dz^1+\cos(t) dz^2\\
& \widetilde \alpha^3=dt,
\end{split}
\end{equation*}
where $t\in [0,2\pi]$ is the fiber coordinate over $z\in \Lambda$.

\quad {\it Example 7.2.}

\quad A more general example is the local form of a 
metric $g=u^2\widetilde g$ conformal 
to the flat case discussed 
above, where $u$ is a smooth function on $\Lambda$. In this case we have 
$g=(\eta^1)^2+(\eta^2)^2$, where 
\begin{equation*}
\eta^1=u \textrm{d}z^1,\qquad  \eta^2=u \textrm{d}z^2.
\end{equation*}
\quad By exterior differentiation it follows
\begin{equation*}
\begin{split}
& a=-\frac{1}{u^2}\frac{\partial u}{\partial z^2}\\
& b=\ \frac{1}{u^2}\frac{\partial u}{\partial z^1}.
\end{split}
\end{equation*}
\quad If we denote by
$\nu:{\mathcal F}_{\textrm{on}}(\Lambda)\to \Lambda$
the bundle of $g$-oriented orthonormal frames on $\Lambda$, we obtain on 
${\mathcal F}_{\textrm{on}}(\Lambda)$ the tautological 1-forms
\begin{equation*}
\begin{split}
& \alpha^1=u \widetilde \alpha^1\\
& \alpha^2=u \widetilde \alpha^2\\
& \alpha^3= \widetilde \alpha^3-*d(\log u),
\end{split}
\end{equation*}
where $*$ is the Hodge operator,  $\widetilde \alpha^1, \widetilde \alpha^2$ and $\widetilde \alpha^3$ are the the tautological 
1-forms and the Levi-Civita connection form of the flat metric $\widetilde g$, respectively.\\
\quad A straightforward computation shows that the Gauss curvature $R$ of $g$ is
\begin{equation}\label{Gauss_curv}
R=-\frac{1}{u^2}\Delta(\log u),
\end{equation}
where $\Delta$ is the Laplace operator in the coordinates $(z^1,z^2)$.\\
\quad It follows that a local form for the coframe $( \alpha^1,  \alpha^2,  \alpha^3)$ is given by
\begin{equation*}
\begin{split}
&  \alpha^1=u\Bigl(\cos(t) dz^1-\sin(t) dz^2\Bigr)\\
&  \alpha^2=u\Bigl(\sin(t) dz^1+\cos(t) dz^2\Bigr)\\
&  \alpha^3=dt-*d(\log u),
\end{split}
\end{equation*}
where $t\in [0,2\pi]$ is the fiber coordinate over $z\in \Lambda$. Here, we denote the 
pullback $\nu^*(u)$ of $u$ to $\mathcal{F}(\Lambda)$ by the same letter.

\subsection{The frame bundle $\mathcal{F}(\Lambda)$}
\quad We return to our setting in \S 7.1, and start with an arbitrary 
Riemannian surface $(\Lambda,g)$ with the area 2-form $dA$ given such that
\begin{equation*}
g=(\eta^1)^2+(\eta^2)^2,\qquad dA=\eta^1\wedge \eta^2>0,
\end{equation*}
where $\{\eta^1,\eta^2\}$ is an $g$-orthonormal coframe on $\Lambda$, and 
$\{e_1,e_2\}$ is its dual frame.\\
\quad We construct as above the $g$-oriented frame bundle $\nu:\mathcal F(\Lambda)\to \Lambda$, where $(z,e_{1\ |
z},e_{2\ |z})$ is a $g$-oriented frame on $\Lambda$.\
\quad Let us denote by $\hat{l}$ the mapping 
\begin{equation*}
\hat{l}:\Sigma\to\mathcal F(\Lambda),\quad
u\mapsto \hat{l}(u)=\Bigl(l(u);l_{*,u}(f_{1\ |u}), l_{*,u}(f_{2\ |u})  \Bigr),
\end{equation*}
where $f_1$, $f_2$ are given in (\ref{frame_change}).

{\bf Proposition 7.1.} {\it The mapping 
$\hat{l}:\Sigma\to\mathcal F(\Lambda)$
defined above is a local diffeomorphism.}

\quad We will give the proof of this result below.\\
\quad We have therefore the commutative diagram. 
\begin{equation*}
\begin{matrix}
\Sigma &  \xrightarrow{\hat{l}}       & \mathcal F(\Lambda)\\
       & l \searrow & \downarrow \nu \\
       &            &  \Lambda
\end{matrix}
\end{equation*}
\quad Remark that due to Proposition 7.1 we can locally identify $\Sigma$ with $\mathcal F(\Lambda)$ as well as the 
coframes $\theta$ and $\alpha$. In order to avoid confusion we will still write $\hat{l}^*$, but we will consider all the 
formulas proved above for the coframe $\theta$ to hold good for $\alpha$ as well via  $\hat{l}^*$.\\
\quad Let us consider now the tautological 1-forms $\{\alpha^1,\alpha^2\}$ on
$\mathcal F(\Lambda)$, i.e.
\begin{equation*}
\nu^*(\eta^1)=\alpha^1,\qquad \nu^*(\eta^2)=\alpha^2,
\end{equation*}
or, equivalently,
\begin{equation}\label{alpha12_omega}
\hat{l}^*(\alpha^1)=m\omega^2,\qquad \hat{l}^*(\alpha^2)=\omega^3.
\end{equation}
\quad A simple computation shows that we must also have
\begin{equation*}
{l}^*(\eta^1)=m\omega^2,\qquad {l}^*(\eta^2)=\omega^3.
\end{equation*}
\subsection{ The structure equations }

\quad We are going to discuss the structure equations on $\mathcal F(\Lambda)$ and $\Lambda$, respectively.

\boxed{
{\it "Upstairs"}
}

\quad We have mentioned already the {\it "upstairs"} structure equations
 on $\mathcal F(\Lambda)$.
%
\quad If we pullback the first two equations to $\Sigma$ by the means of $\hat{l}^*$, it follows
\begin{equation*}
\begin{split}
&d (\hat{l}^* \alpha^1)=\hat{l}^*(\alpha^2)\wedge\hat{l}^*(\alpha^3)\\
&d (\hat{l}^* \alpha^2)=\hat{l}^*(\alpha^3)\wedge\hat{l}^*(\alpha^1)
\end{split}
\end{equation*}
and from here, by using (\ref{alpha12_omega}) we obtain
\begin{equation}\label{alpha3_omega}
\hat{l}^*(\alpha^3)=m\omega^1+m_3\omega^2
\end{equation}
on $\Sigma$.\\
\quad Remark that 
\begin{equation*}
\hat{l}^*(\alpha^1\wedge\alpha^2\wedge  \alpha^3)=
m^2\omega^1\wedge\omega^2\wedge\omega^3\neq 0,
\end{equation*}
i.e. $\hat{l}$ is indeed a local diffeomorphism and this proves the Proposition 7.1 above.

\boxed{
{\it "Downstairs"}
}

\qquad The {\it "downstairs"} structure equations on $\Lambda$ are
\begin{equation*}
\begin{split}
& d\eta^1=\eta^2\wedge \eta^3\\
& d\eta^2=\eta^3\wedge \eta^1\\
& d\eta^3=R\ \eta^1\wedge \eta^2,
\end{split}
\end{equation*}
where $R$ is the {\it "downstairs"} Gauss curvature of $(\Lambda,g)$.\\
\quad We pullback the last equation above to $\mathcal F(\Lambda)$ by means of $\nu^*$. It follows
\begin{equation*}
d\alpha^3=\nu^*(R\ \eta^1\wedge \eta^2).
\end{equation*}
\quad On the other hand, by exterior differentiation of (\ref{alpha3_omega}) we obtain
\begin{equation*}
\begin{split}
\hat{l}^*(d\alpha^3) & =d(m\omega^1)+d(m_3\omega^2)=(m-m_{33})\omega^2\wedge\omega^3\\
&=\frac{m-m_{33}}{m}\hat{l}^*(\alpha^1)\wedge\hat{l}^*(\alpha^2)=
(1-\frac{m_{33}}{m})l^*(\eta^1\wedge\eta^2).
\end{split}
\end{equation*}
\quad It follows
\begin{equation*}
l^*(R\eta^1\wedge\eta^2)=(1-\frac{m_{33}}{m})l^*(\eta^1\wedge\eta^2),
\end{equation*}
and from here we obtain the following {\it curvature condition} on $\Sigma$:
\begin{equation}\label{condC_up}
(C)\qquad \frac{m_{33}}{m}l^*(\eta^1\wedge\eta^2)=l^*\Bigl[(1-R)\eta^1\wedge\eta^2\Bigr].
\end{equation}
\quad We would like to express now the quantity $\frac{m_{33}}{m}$ living on $\Sigma$ as the image of a quantity living 
on $\Lambda$ through $l^*$. \\
\quad Recall from the general theory that if $\{e_1,e_2\}$ is an adapted frame to the geometry of the Riemannian 
surface $(\Lambda,g)$, this is equivalent with giving a section of the frame bundle $\nu:\mathcal F(\Lambda)\to \Lambda
$, i.e.
\begin{equation*}
s:\Lambda\to \mathcal F(\Lambda),\qquad \nu\circ s=id_\Lambda,
\end{equation*}
i.e. we have a so called {\it first order adapted lift}.\\
\quad Let us consider next an arbitrary smooth function $\bar{m}$ on $\Lambda$, and lift it {\it ``upstairs''}, i.e. we obtain 
a function $\widetilde{m}=\bar{m}\circ \nu$ on $\mathcal F(\Lambda)$, such that $s^*(\widetilde{m})=\bar{m}$, and a function $m$ on $
\Sigma$ such that
\begin{equation}
m=\hat{l}^*(\widetilde{m})=\hat{l}^*(\nu^*\bar{m})=(\nu\circ \hat{l})^*\bar{m}.
\end{equation}
\quad We take next the exterior derivative of the relation $m=l^*(\bar m)$. It follows
\begin{equation*}
\begin{split}
dm=l^*(d\bar m)=l^*(\bar m_1\eta^1+\bar m_2\eta^2)=l^*(\bar m_1)m\omega^2+l^*(\bar m_2)\omega^3,
\end{split}
\end{equation*}
i.e. $dm$ is a linear combination of the 1-forms $\omega^2$, $\omega^3$. This implies
\begin{equation*}
m_1=0.
\end{equation*}
\quad It follows that this $m$ can be used to relate the coframes $\omega$ and $\alpha$ as in \S6.1. Under these
conditions, we take the exterior derivative of the relation 
$m=\hat{l}^*(\widetilde{m})$. It follows that
\begin{equation*}
\begin{split}
dm & =m_2\omega^2+m_3\omega^3=\hat{l}^*(\widetilde{m}_1\alpha^1+\widetilde{m}_2\alpha^2+
\widetilde{m}_3\alpha^3)\\
& =\hat{l}^*(\widetilde{m}_1)m\omega^2+\hat{l}^*(\widetilde{m}_2)\omega^3+\hat{l}^*(\widetilde{m}_3)
(m \omega^1+m_3\omega^2),
\end{split}
\end{equation*}
and from here, we obtain
\begin{equation*}
\begin{split}
& \hat{l}^*(\widetilde{m}_1)=\frac{m_2}{m}\\
& \hat{l}^*(\widetilde{m}_2)=m_3\\
& \hat{l}^*(\widetilde{m}_3)=0.
\end{split}
\end{equation*}
\quad Remark that Proposition 7.1 together with the last condition above imply that 
 \begin{equation*}
\widetilde{m}_3=0.
\end{equation*}
\quad By a straightforward computation we also obtain
 \begin{equation*}
\hat{l}^*(\widetilde{m}_{22})=m_{33}.
\end{equation*}
\quad Recall that $(\eta^1,\eta^2,\eta^3)=s^*(\alpha^1,\alpha^2,\alpha^3)$, and using now the relation $
\bar{m}=s^*(\widetilde{m})$ we have
\begin{equation*}
s^*(d\widetilde{m})=s^*(\widetilde{m}_1)\eta^1+s^*(\widetilde{m}_2)\eta^2,
\end{equation*}
where we have put $d\widetilde{m}=\widetilde{m}_1\alpha^1+\widetilde{m}_2\alpha^2$ on $\mathcal F(\Lambda)$ and 
$d\bar m=\bar m_1\eta^1+\bar m_2\eta^2$ on $\Lambda$.\\
\quad Then, it follows
\begin{equation*}
\begin{split}
& \bar m_1=s^*(\widetilde{m}_1)\\
& \bar m_2=s^*(\widetilde{m}_2).
\end{split}
\end{equation*}
%
\quad A straightforward computation using (\ref{lands_cond_theta}), (\ref{Ricci_cond_theta}) pulled back through $\hat{l}^*$ shows 
that
\begin{equation*}
d\widetilde{m}_2=\widetilde{m}_{22}\alpha^2+\widetilde{m}_1\alpha^3,
\end{equation*}
and pulling this equation back through $s^*$ we get
\begin{equation*}
s^*(\widetilde{m}_{22})=\bar m_{22}+b\bar m_1,
\end{equation*}
where $b$ is the function  on $\Lambda$ from $d\eta^2=b\eta^1\wedge\eta^2$.\\
\quad In the same way we obtain
\begin{equation*}
\begin{split}
&  s^*(\widetilde{m}_{11})=\bar m_{11}-a\bar m_2,\\
&  s^*(\widetilde{m}_{12})= s^*(\widetilde{m}_{21})=\bar m_{12}-b\bar m_2=\bar m_{21}+a\bar m_1,
\end{split}
\end{equation*}
where we take into account the Ricci type identity on $\Lambda$:
\begin{equation*}
\bar{m}_{21}-\bar{m}_{12}+a\bar m_1+b\bar m_2=0.
\end{equation*}
\quad Hence, we obtain
\begin{equation*}
m_{33}=l^*(s^*(\widetilde{m}_{22}))=l^*(\bar m_{22}+b\bar m_1).
\end{equation*}
\quad Using now this in (\ref{condC_up}) we are led to the following {\it curvature relation on $\Lambda$}:
\begin{equation}\label{condC_down}
(C) \qquad \frac{\bar m_{22}+b\bar m_1}{\bar m}=1-R,
\end{equation}
which, together with the {\it Landsberg condition on $\Lambda$} , namely
\begin{equation}\label{condL_down}
(L) \qquad \bar m_{12}-b\bar m_2=\bar m_{21}+a\bar m_1=0,
\end{equation}
are the fundamental relations to be satisfied by $\bar m$ on $\Lambda$.\\
\quad Remark that the non-triviality relations $m_2\neq 0$, $m_3\neq 0$ are equivalent to 
\begin{equation*}
\widetilde{m}_1\neq 0,\qquad \widetilde{m}_2\neq 0
\end{equation*}
on $\mathcal F(\Lambda)$
or, equivalently,
\begin{equation}\label{condN_down}
(N)\qquad \bar{m}_1\neq 0,\qquad \bar{m}_2\neq 0
\end{equation}
on $\Lambda$.

\section{Constructing local generalized unicorns}
\subsection{Recovering the generalized Landsberg structure}
\quad Conversely, one can locally construct a generalized Landsberg structure as follows. Let us consider
\begin{enumerate}
\item an oriented Riemannian surface $(\Lambda,g)$ of Gauss curvature $R$, and
\item a function $\bar m$ on $\Lambda$ that satisfies the PDE system (\ref{condC_down}), (\ref{condL_down})
 with the non-triviality conditions (\ref{condN_down}).
\end{enumerate}
\quad Then, on the orthonormal frame bundle $\nu:\mathcal F(\Lambda)\to \Lambda$ there exist the tautological 1-forms 
$\alpha^1$, $\alpha^2$ and the Levi-Civita connection form $\alpha^3$ that satisfy the usual structure equations
\begin{equation}\label{riemann_struct_eq}
\begin{split}
d\alpha^1& =\alpha^2\wedge\alpha^3\\
d\alpha^2& =\alpha^3\wedge\alpha^1\\
d\alpha^3& =\nu^*(R)\ \alpha^1\wedge\alpha^2.
\end{split}
\end{equation}
\quad Let us construct the coframing 
\begin{equation}\label{inverse_coframe}
\begin{split}
& \bar{\omega}^1=\frac{1}{\widetilde{m}}(\alpha^3-\frac{\widetilde{m}_2}{\widetilde{m}}\alpha^1)\\
& \bar{\omega}^2=\frac{1}{\widetilde{m}}\alpha^1\\
& \bar{\omega}^3=\alpha^2,
\end{split}
\end{equation}
where $\widetilde{m}=\nu^*(\bar m)$.\\
\quad It follows from Section 6, Section 7 that $\{\bar\omega^1,\bar\omega^2,\bar\omega^3\}$ is  a non-trivial 
generalized Landsberg structure on the 3-manifold $\mathcal F(\Lambda)$ with the invariants
\begin{equation*}
I=-2\frac{\widetilde{m}_2}{\widetilde{m}},\qquad K=\widetilde{m}^2.
\end{equation*}
\quad By similar computations as in Section 4 one can show by means of 
Cartan-K\"ahler theorem that the PDE system (\ref{condC_down}), (\ref{condL_down}) is involutive. 
We will not discuss here the most general situation, but a particular case will be described below.
We recall also that a Riemannian structure on a surface depends on a function of two variables, say $u$ on $\Lambda$ (this is a consequence of the existence of isothermal coordinates on a Riemannian surface).\\
\quad Summarizing, it follows from the Cartan--K\"ahler theorem used in Section 4 that the degree of freedom of the scalar invariants 
$I$, $K$ of a generalized Landsberg 
structure locally depends on two arbitrary functions of two variables (see \S 4.1, \S 4.2). 
We point out that these two functions of two variables are in the Cartan-K\"ahler sense, i.e. they show the degree of freedom of $(I,K)$, but one should not think that they are exactly 
the functions $u$ and $\bar{m}$ used in the precedent section.  \\
\quad More generally, a generalized Landsberg structure, i.e. the 
coframe $\{\omega^1,\omega^2,\omega^3\}$ together with the scalar invariants $I$, $K$, depends on 3 functions of 3 variables (see \S 4.3).
 A particular case is the generalized Landsberg 
structure (\ref{inverse_coframe}) constructed using a function $u$ on $\Lambda$, from the Riemannian structure $(\Lambda,g)$ downstairs, and a function $\bar m$ on $\Lambda$ satisfying  (\ref{condC_down}), (\ref{condL_down}). We will show in the next section that the degree of freedom of the pair of 
functions $(u,\bar m)$ is actually 4 functions of 1 variable (see Proposition 8.1).\\
\quad Remark that our solution has a lower degree of freedom than the general solution predicted by our first use of 
Cartan--K\"ahler theorem in Section 4 due to our particular choice of the coframe changing (\ref{inverse_coframe}), so there is 
no contradiction with our results in Section 4.\\
\quad Remark also that our condition $\widetilde m_1=0$ implies that the directional derivative of the invariant $K$ with 
respect to $\widetilde \omega^1$ vanishes, in other words we are considering here an integral manifold of the linear Pfaffian 
system
(\ref{Lands_Pfaffian}) passing through the initial condition 
\begin{equation*}
(u_0,I(u_0),K(u_0),I_1(u_0),I_3(u_0),0,K_2(u_0)),
\end{equation*}
as explained in \S 4.2, where the invariants $I$, $K$ are given above.

\subsection{A local form}
\quad In order to construct a local form for the   generalized Landsberg structure given by (\ref{inverse_coframe}), 
we 
are going to use Zoll projective structures.\\
\quad Let us start with a Riemannian metric $g=u^2[(dz^1)^2+(dz^2)^2]$ on the surface $\Lambda$ 
with the 
Christoffel symbols $\Gamma_{jk}^i$, and construct the 1-form $\gamma$ on $\Lambda$ as in \eqref{Gamma}, 
\eqref{small_gamma}.\\
\quad By putting
  $\gamma=d(\log u)$, i.e.  
\begin{equation}\label{u_ODE}
\frac{1}{u}\frac{\partial u}{\partial z^i}=\gamma_i,\qquad i=1,2,
\end{equation}
in some isothermal coordinates $(z^1,z^2)\in \Lambda$, it follows that the Gauss curvature of the Riemannian metric 
$g=u^2[(dz^1)^2+(dz^2)^2]$ will be given by 
\begin{equation*}
R=-\frac{1}{u^2}\textrm{div}\gamma
\end{equation*}
as explained in \S 3.2. See also Example 7.2 for other formulas.\\
\quad On the other hand, in order to obtain a generalized Landsberg structure upstairs, we need a function $\bar{m}$ on 
$\Lambda$ that satisfies the conditions (\ref{condC_down}), (\ref{condL_down})
 and the  non-triviality conditions (\ref{condN_down}). \\
\quad If we denote by numerical subscripts the directional derivatives of $\bar{m}$ with respect to the $g$-orthonormal 
coframe 
\begin{equation*}
\eta^1=udz^1,\quad \eta^2=udz^2,
\end{equation*}
and with letters the partial derivatives, then straightforward computations show the expression of first order directional 
derivatives
\begin{equation}\label{m_1st_deriv}
\bar{m}_i=\frac{1}{u}\bar{m}_{z^i},\qquad i=1,2,
\end{equation}
and second order directional derivatives
\begin{equation}\label{m_2st_deriv}
\begin{split}
& \bar{m}_{11}=\frac{1}{u^2}(-\gamma_1\bar{m}_{z^1}+\bar{m}_{z^1z^1})\quad
 \bar{m}_{12}=\frac{1}{u^2}(-\gamma_2\bar{m}_{z^1}+\bar{m}_{z^1z^2})\\
& \bar{m}_{21}=\frac{1}{u^2}(-\gamma_1\bar{m}_{z^2}+\bar{m}_{z^2z^1})\quad
 \bar{m}_{22}=\frac{1}{u^2}(-\gamma_2\bar{m}_{z^2}+\bar{m}_{z^2z^2}).
\end{split}
\end{equation}
\quad It follows from (\ref{condC_down}), (\ref{condL_down}) that $\bar{m}$ must satisfy
\begin{enumerate}
\item {\it The Landsberg condition}
\begin{equation}\label{condL_2}
(L)\qquad \bar{m}_{z^1z^2}=\gamma_1\bar{m}_{z^2}+\gamma_2\bar{m}_{z^1},
\end{equation}			
\item {\it The curvature condition}
\begin{equation}\label{condC_2}
(C)\qquad \bar{m}_{z^2z^2}= -(\gamma_1m_{z1}-   \gamma_2m_{z2})+u^2+\textrm{div}\gamma.
\end{equation}
\end{enumerate}
\quad It follows that these two conditions can be regarded as a PDE system for $\bar{m}$ on $\Lambda$, 
where $\gamma$'s is given by (\ref{u_ODE}).\\
\quad The first question that arises is the involutivity of such a PDE system. We will discuss this using our favorite tool, 
the Cartan-K\"ahler theorem. \\
\quad Let $J^2(\mathbb{R}^2, \mathbb{R}^2)$ be a second order jet space of 
two functions on a plane.
The second jet space $J^2 (\mathbb{R}^2, \mathbb{R}^2)$ has the canonical system 
\begin{equation*}
C^2= \{\theta_{i}^{j}=0 \quad (i=0,1,2 ,j=1,2)\}
\end{equation*}
 where
$(z^1,z^2,\bar{m},u, \bar{m}_{z^1},\bar{m}_{z^2},
u_{z^1},u_{z^2},\bar{m}_{z^1 z^1},\bar{m}_{z^1 z^2},\bar{m}_{z^2 z^2},
u_{z^1 z^1},u_{z^1z^2},u_{z^2 z^2})$ are the coordinates on $J^2(\mathbb{R}^2, \mathbb{R}^2)$ 
and
\begin{eqnarray*}
\theta_{0}^{1}=d\bar{m}-\bar{m}_{z^1}dz^1 -\bar{m}_{z^2}dz^2  \qquad
&,&\quad  \theta_{0}^{2}= du-u_{z^1}dz^1 -u_{z^2}dz^2\ , \\
\theta_{1}^{1}=d\bar{m}_{z^1}-\bar{m}_{z^1 z^1}dz^1 -\bar{m}_{z^1 z^2}dz^2  
  &,& \quad\theta_{1}^{2}= du_{z^1}-u_{z^1 z^1}dz^1 -u_{z^1 z^2}dz^2\ , \\
\theta_{2}^{1}= d\bar{m}_{z^2}-\bar{m}_{z^1 z^2}dz^1 -\bar{m}_{z^2 z^2}dz^2 &,& \quad\theta_{2}^{2}= du_{z^2}-u_{z^1 
z^2}dz^1 -u_{z^2 z^2}dz^2
\end{eqnarray*}
are the canonical contact forms. \\
 \quad We consider the system of PDE formed by the equations $(L),(C)$, namely, 
\[
R=\{(L),(C) \} \subset J^2(\mathbb{R}^2, \mathbb{R}^2),\quad I= C^2 |_R , \quad
 \Omega =dz^1 \wedge dz^2,
\]
with coordinates $(z^1,z^2,\bar{m},u, \bar{m}_{z^1},\bar{m}_{z^2},
u_{z^1},u_{z^2},\bar{m}_{z^1 z^1},
u_{z^1 z^1},u_{z^1z^2},u_{z^2 z^2})$ on $R$.\\
\quad By a straightforward computation we find that 
the Pfaffian system $I$ has absorbable torsion. Moreover, its tableau is given by
\begin{equation}
\begin{pmatrix}
0 & \qquad & 0 \\
a & \qquad & 0 \\
0 & \qquad & \frac{\bar m}{u}(b+d) \\
0 & \qquad & 0 \\
b & \qquad & c \\
c & \qquad & d
\end{pmatrix}
\end{equation}
and the characters of the tableau are $s_1=4,\ s_2=0$. Since  the dimension of 
the space of integral elements is $4=s_1+2s_2$, Cartan's Test for 
involutivity implies that the system is involutive. \\
\quad Hence, in the analytic category, the Cartan-K\"ahler theorem implies that 
the solutions exist, and, roughly speaking, they depend on 4 functions of 1 variable. \\
\quad We are led in this way to the following result.

{\bf Proposition 8.1.} \\
{\it \quad The system of partial differential equations (L), (C) for two unknown functions u, $\bar m$ of two variables has solutions.
Moreover, these solutions depend in Cartan-K\"ahler sense on 4 functions of 1 variable.}

\quad We obtain therefore the following prescription for constructing generalized Landsberg structures:

\quad $\bullet$ Start with a smooth surface $\Lambda$ with local coordinates $z^1$, $z^2$ and consider the functions $
\bar{m}, u:\Lambda\to \mathbb R$  which satisfy (\ref{condL_2}),  (\ref{condC_2}). The existence of such an $\bar{m}$ and $u$ is 
guaranteed by the Cartan-K\"ahler theorem (Proposition 8.1).

\quad $\bullet$ Denote by $g=u^2[(dz^1)^2+(dz^2)^2]$ the corresponding Riemannian metric on $\Lambda$ conformal 
equivalent to the flat metric, and by $R$ its Gauss curvature given by (\ref{Gauss_curv});

\quad $\bullet$ Construct the $g$-orthonormal frame bundle $\nu:\mathcal F(\Lambda)\to \Lambda$ with the tautological 
1-forms $\alpha^1$, $\alpha^2$ and the Levi-Civita connection form $\alpha^3$;

\quad $\bullet$ Lift the function $\bar{m}$ to $\Sigma:=\mathcal F(\Lambda)$ as
$\widetilde m:=\nu^*(\bar{m})$;

\quad $\bullet$ Construct the coframe $(\bar{\omega}^1,\bar{\omega}^2,\bar{\omega}^3)$ on $\Sigma=\mathcal 
F(\Lambda)$ given by (\ref{inverse_coframe}).

\quad Then, we have

\quad {\bf Theorem 8.2.} {\it The coframe $(\bar{\omega}^1,\bar{\omega}^2,\bar{\omega}^3)$ constructed above is a 
generalized Landsberg structure on the 3-manifold $\Sigma=\mathcal F(\Lambda)$.}

\quad Indeed, remark first that $\widetilde m:=\nu^*(\bar{m})$ implies $s^*(\widetilde m)=\bar{m}$, as well as $\widetilde m_3=0$ by taking the exterior 
derivative. Then, in the present setting, similar computations with those in \S 7.2 show that conditions (L) and (C) upstairs in Proposition 6.2 hold good. Computing 
now the structure equations of the coframe $\bar{\omega}$ and making use of (\ref{riemann_struct_eq}) 
and properties in Proposition 6.2, one can easily verify that $\bar{\omega}$ is a generalized Landsberg structure on the 3-manifold 
$\Sigma=\mathcal{F}(\Lambda)$.\\
\quad Using the normal form from Example 7.2 in \S 7.1, we obtain the following normal form of this generalized unicorn:
\begin{equation}\label{normal_form}
\begin{split}
& \bar{\omega}^1=\frac{1}{\widetilde{m}}\Bigl[dt-*d(\log u)-\frac{u\ \widetilde{m}_2}{\widetilde{m}}\Bigl(\cos(t) 
dz^1-\sin(t) dz^2\Bigr)\Bigr]\\
& \bar{\omega}^2=\frac{u}{\widetilde{m}}\Bigl(\cos(t) dz^1-\sin(t) dz^2\Bigr)\\
& \bar{\omega}^3=u\Bigl(\sin(t)dz^1+\cos(t)dz^2\Bigr),
\end{split}
\end{equation}
where $\widetilde{m}=\nu^*(\bar m)$, 
$\widetilde{m}_2=\nu^*(\frac{1}{u}\frac{\partial \bar{m}}{\partial z^2})$
and 
$t\in [0,2\pi]$ is the fiber coordinate over $z=(z^1,z^2)\in \Lambda$. Here, we denote again the prolongation $\nu^*(u)$ of $u$ to $\mathcal{F}(\Lambda)$ by the same letter.

\section{Concluding remarks}
\quad In the present note we have shown how is possible to construct a non-trivial generalized Landsberg structure 
$\{\omega^1,\omega^2,\omega^3\}$ on a 3-manifold $\Sigma$ using a Riemannian metric $g$ on a surface $\Lambda$ 
that basically depends on 2 functions of 1 variable, namely, $u$ and $\bar{m}$. Due to Cartan-K\"{a}hler Theorem in \S 8.2, we 
know that these functions are locally described by 4 functions of 1 variable, case included in 
the general solution predicted by 
Cartan-K\"ahler Theory in Section 4. A local form of it is given by (\ref{normal_form}). This generalized Landsberg structure is 
{\it locally amenable} in the sense of \S 5.2. Our generalized unicorn has the fundamental geometrical property that 
its indicatrix foliation $\{\omega^1=0,\omega^2=0\}$ coincides with the  geodesic foliation $\{\alpha^1=0,\alpha^3=0\}$ of 
the Riemannian metric $g$ of $\Lambda$.\\
\quad However, our initial intention was to search for classical unicorns on surfaces, i.e. generalized Landsberg 
structures that satisfy the conditions of Theorem 2.1.\\
\quad Recall that a generalized Finsler structure is amenable if the indicatrix foliation 
$\mathcal Q=\{\omega^1=0,\omega^2=0\}$ is amenable, i.e. the leaf space is a differentiable manifold. \\
\quad Let us also recall that a Zoll metric on $S^2$ depends on one odd arbitrary function on one variable (see 
\cite{B1978} and \cite{LM2002} for details). We are lead in this way to the following 

\quad {\bf Conjecture 9.1.} {\it There exists a solution u of (\ref{condL_2}), (\ref{condC_2}) that gives a Riemannian metric 
$g=u^2[(dz^1)^2+(dz^2)^2]$ whose Levi-Civita connection $\nabla^g$ belongs to a Zoll projective class on $S^2$.
}

\quad If we accept this conjecture as true, then we just have constructed a generalized Landsberg structure $\{\bar 
\omega^1,\bar \omega^2,\bar \omega^3\}$ on the frame bundle $\Sigma:=\mathcal F(S^2)$ of a Riemannian surface $
(S^2,g)$ whose Levi-Civita connection $\nabla^g$ belongs to a Zoll projective structure on $S^2$, in other words, the
geodesic foliation $\mathcal P=\{\alpha^1=0,\alpha^3=0\}$ of $g$ foliates the 3-manifold $\Sigma$ by circles. Remark in 
the same time that we had constructed our coframe $\bar \omega$ from 
$\alpha$ by (\ref{inverse_coframe}) such that  its indicatrix foliation $\mathcal Q=\{\omega^1=0,\omega^2=0\}$ coincides 
with the geodesic foliation $\mathcal P=\{\alpha^1=0,\alpha^3=0\}$ of $g$. Then, by the properties of Zoll projective 
structure on $S^2$ described partially in \S 3.2 it follows that the space of geodesics, say $M$, of the metric $
(\Lambda=S^2,g)$ is a differentiable manifold, and hence, the generalized Landsberg structure  $\{\bar \omega^1,\bar 
\omega^2,\bar \omega^3\}$ is globally amenable. In other words, the map $\pi:\Sigma\to S^2$ is a smooth submersion. 
Obviously, the leaves of the indicatrix foliation  $\{\bar \omega^1=0,\bar \omega^2=0\}$ are diffeomorphic to $S^1$, so 
they must be compact.\\
\quad Finally, in order to have a true classical unicorn, we have to show more, namely that the canonical immersion $
\iota:\Sigma\to TM$, given by $\iota(u)=\pi_{*,u}(\hat{e}_2)$ is injective on each $\pi$-fiber $\Sigma_x$, as stated in 
Theorem 2.1. This is not so difficult to prove. Let us denote by 
\begin{equation*}
\gamma_u:[a,b]\to \Sigma
 \end{equation*}
the geodesic flow of the Zoll projective structure $[\nabla]$ on $S^2$  through the point $u\in \Sigma$, and let us take 
another point, say $u_1$ on the same leaf, i.e. there exist some parameter values $s_0, s_1\in [a,b]$ such that
\begin{equation*}
\gamma_u(s_0)=u,\qquad \gamma_u(s_1)=u_1
\end{equation*}
on $\Sigma$. \\
\quad From \S 3.2 we know that the leaves $\gamma$ are closed, periodic, simple curves of same length on $
\Sigma$, i.e. for 
\begin{equation*}
\gamma_u(s_0)=u \neq \gamma_u(s_1)=u_1 \Longrightarrow \hat{e}_{2\ |\gamma_u(s_0)}\neq \hat{e}_{2\ |
\gamma_u(s_1)},
\end{equation*}
where $\hat{e}_{2}\in T_{\gamma}\Sigma$ is thought as a vector field along $\gamma$. Applying to this the linear map $
\pi_{*,u}$ it follows
\begin{equation*}
\pi_{*,\gamma(s_0)}(\hat{e}_{2\ |\gamma_u(s_0)})\neq \pi_{*,\gamma(s_1}(\hat{e}_{2\ |\gamma_u(s_1)})
\end{equation*}
and therefore it follows that $\iota$ must be injective on each $\pi$-fiber $\Sigma_x$.\\
\quad Then, from Theorem 2.1 we can conclude

{\it \quad There are Landsberg structures on $M=S^2$ which are not Berwald type, provided the conjecture above is 
true.}


\section{Appendix. The Cartan--K\"ahler theorem for linear Pfaffian systems}

\quad We give a short outline of the main tool used in the present paper, the Cartan--K\"ahler theorem for linear Pfaffian 
systems. This theorem is presented in several textbooks, for \cite{Br et al 1991}, \cite{IL2003}, \cite{O1995}, etc., but our 
presentation here follows our favorite monograph \cite{IL2003}.\\
\quad Let us denote by $\Omega^*(\Sigma)=\bigoplus_{k} \Omega^k(\Sigma)$ the space of smooth differential forms on 
the manifold $\Sigma$. It is a standard fact that $\Omega^*(\Sigma)$ is a graded algebra under the wedge product.  \\
\quad A subspace $\mathcal{I}\subset \Omega^*(\Sigma)$ is called {\it an exterior ideal} or {\it an algebraic ideal} if it is a 
direct sum of homogeneous subspaces
(namely, $\mathcal{I} = \bigoplus_{k} {\mathcal{I}}^k$, 
${\mathcal{I}}^k \subset \Omega^k(\Sigma)$.) and it satisfies 
\[
\omega\wedge \eta\in \mathcal{I},
\]
 for $\omega\in \mathcal{I}$ and {\it any} differential form $\eta\in \Omega^*(\Sigma)$.\\
\quad An exterior ideal is called {\it a differential ideal} if for any $\omega\in \mathcal{I}$, we have $d\omega\in 
\mathcal{I}$ also. \\
\quad A differential ideal $\mathcal{I} \subset \Omega^*(\Sigma)$ is called an 
{\it exterior differential system} on a manifold $\Sigma$ (EDS for short).\\
\quad A set of differential forms of arbitrary degree $\{\omega^1,
\omega^2, \dots,\omega^k \}$ is said to 
{\it generate the EDS} $\mathcal{I}$ if any $\theta \in \mathcal{I}$ 
can be written as a finite ``linear combination'', namely 
\begin{equation*}
\mathcal{I} =\{ \sum_{i=1}^{k} \alpha^i \wedge \omega^i+
\sum_{i=1}^{k} \beta^i \wedge d\omega^i \ |\ \alpha^i,\beta^i 
\in \Omega^*(\Sigma)
\}.
\end{equation*}
\quad A {\it Pfaffian system} $\mathcal{I}$ on a manifold $\Sigma$ is an EDS 
finitely generated by 1-forms $\{\omega^1,\omega^2,\dots\omega^k\}$ only.\\
\quad For an EDS $\mathcal{I}$ on a manifold $\Sigma$, a decomposable differential $k$-form 
$\Omega$ (up to scale) 
is called the independence condition if $\Omega$ does not vanish
modulo $\mathcal{I}$ on $\Sigma$.\\
\quad We denote by $(\mathcal{I}, \Omega)$ a pair of an EDS and an 
independence condition on a manifold $\Sigma$.\\
\quad A submanifold $f:M\to\Sigma$ is called  {\it an integral submanifold} (or solution) of the EDS $(\mathcal{I}, \Omega)
$ if
\begin{equation*}
\begin{split}
& f^*(\theta^a)=0, \qquad \theta^a\in \mathcal{I},\\
& f^*(\Omega)\neq 0.
\end{split} 
\end{equation*} 
\quad Remark also that $f^*(\theta)=0$ imply $f^*(d\theta)=0$.\\
\quad There is a notion of infinitesimal solution also. A $k$-dimensional subspace $E\subset T_x\Sigma$ is called {\it an 
integral element} of 
$(\mathcal{I}, \Omega)$ if 
\begin{equation*}
\begin{split}
& \theta^a_{|_E}=0, \qquad \theta^a\in \mathcal{I},\\
& \Omega_{|_E}\neq 0.
\end{split} 
\end{equation*} 
\quad Usually one regards $E$ as an element of the Grassmannian $G_k(T_x\Sigma)$ of $k$-planes through the origin 
of the vector space $T_x\Sigma$. The space of
$k$-dimensional integral elements of $(\mathcal{I},\Omega)$ is usually denoted by $\mathcal{V}_k(\mathcal{I},\Omega)
$.\\
\quad Roughly speaking, a differential system will be called {\it integrable} if one can determine its integral manifolds of a 
prescribed dimension passing through each point. 
In the case of a  Pfaffian system with the maximum degree 
independence condition, its integrability is guaranteed by Frobenius theorem. 
However, in the case when 
the independence condition is not the maximum degree, then one has to use more powerfull  
tools as 
the Cartan-K\"{a}hler Theorem.\\
\quad Let $(I,J)$ be a pair of a collection of $1$-forms $I=\{\theta^1, \theta^2,
\ldots, \theta^s \}$ and $J=\{\omega^1 ,\omega^2 ,\ldots ,\omega^k \}$
which are linearly independent modulo $I$.\\
\quad Remark that $(I,J)$ induces an EDS $(\mathcal{I},\Omega)$ by a Pfaffian system 
$\mathcal{I}$ generated by $I$ and the independence condition $\Omega=
\omega^1 \wedge \omega^2 \wedge \ldots \wedge \omega^k$.\\
\quad The pair $(I,J)$ is called a {\it linear Pfaffian system} if 
\begin{equation*}
 d\theta^a\equiv 0\qquad \mod J,
\end{equation*}
for all $\theta^a$ in $I$.\\
\quad If $(I,J)$ is a linear Pfaffian system, let us denote by 
$\pi^{\epsilon}$, $\epsilon=1,2,\dots,\dim \Sigma-s-k$ such that $T^*\Sigma$ is locally spanned by $\theta^a,\omega^i,
\pi^\epsilon$. The coframing $\theta^a,\omega^i,\pi^\epsilon$ is called {\it adapted} to the filtration 
$I\subset J\subset T^*\Sigma$. 
It follows immediately that there must locally exist some functions $A_{\epsilon i}^a$ and $T_{ij}^a$ on $\Sigma$ such 
that
\begin{equation}\label{struct_eq_torsion}
 d\theta^a\equiv A_{\epsilon i}^a\pi^\epsilon\wedge\omega^i+T_{ij}^a\omega^i\wedge\omega^j\qquad \mod I.
\end{equation}
\quad The terms $T_{ij}^a\omega^i\wedge\omega^j$ in (\ref{struct_eq_torsion}) are called {\it apparent torsion}. Apparent torsion must be normalized before prolonging the system. Namely, one have to choose, if possible, some new one forms $\tilde{\pi}^\epsilon$ such that $\tilde{T}_{ij}^a=0$, with respect to the new coframe $
\theta^a,\omega^i,\tilde{\pi}^\epsilon$ on $\Sigma$. In this case one says that {\it the 
apparent torsion is absorbable}.\\
\quad If this is not possible, then one says that there is {\it torsion} and in this case the system admits no integral elements.\\  
\quad Remark that the functions $A_{\epsilon i}^a$ and $T_{ij}^a$ depend on the choices of the bases for $I$ and $J$. 
However, one can construct invariants from these functions. Indeed, for a fixed generic point $x\in \Sigma$, the {\it 
tableau of $(I,J)$ at x} is defined as
$\Sigma$ such that
\begin{equation*}
 A_x:=\{A_{\epsilon i}^aw_a\otimes v^i\ :\ 1\leq \epsilon\leq \dim \Sigma-\dim J_x\}\subseteq W\otimes V^*,
\end{equation*}
where $V^*:=(J/I)_x$, $W^*=I_x$, $w^a=\theta^a_x$, $v^j=\omega^j_x$. A standard argument of linear algebra shows 
that $A_x$ is independent of any choices.\\
\quad We fix a point $x\in \Sigma$ and denote the tableau $A_x$ simply with $A\in  W\otimes V^*$. The tableau $A$ 
depends on the basis $b=(v^1,v^2,\dots,v^n)$ of $W$. One defines 
\begin{tabular}{rcl}
 $s_1(b)$ & = & no. of independent entries in the first col. of $A$\\
$s_1(b)+s_2(b)$ & = & no. of independent entries in the first 2 col. of $A$\\
\    & \dots & \ \\
$s_1(b)+\dots+s_n(b)$ & = & no. of independent entries in $A$.
\end{tabular}

\quad Equivalently, one can see that {\it the characters $s_1(b),s_2(b),\dots,s_n(b)$ of the tableau A} do not depend 
actually on the choice of the basis $b$ of $W$, but only on the flag of subspaces
\begin{equation*}
 F:\ (0)=F_n\subset F_{n-1}\subset \dots F_1\subset F_0=V^*.
\end{equation*}
\quad This allows us to rewrite $s_k(b)$ as $s_k(F)$. By defining
\begin{equation*}
 A_k(F)=(W\otimes F_k)\cap A,
\end{equation*}
it follows that
\begin{equation*}
 \dim A_k(F)=s_{k+1}(F)+\dots s_n(F).
\end{equation*}
One can easily see that $A_k(F)$ is the subspace of matrices in $A$ for which the first $k$ columns are zero with 
respect to the basis $b$ for $V$.\\
\quad One defines next the {\it reduced characters of the tableau A} as\\
\begin{tabular}{rcl}
 $s_1$ & = & $\max\{s_1(F)\ :\ \text{all flags}\}$\\
$s_2$ & = & $\max\{s_1(F)\ :\ \text{flags with}\  s_1(F)=s_1\}$\\
\    & \dots & \ \\
$s_n$ & = & $\max\{s_n(F)\ :\ \text{flags with} \ s_1(F)=s_1,\dots,s_{n-1}(F)=s_{n-1}\}$.
\end{tabular}

\quad These scalars are invariants of the tableau $A$, i.e. they are independent of any choice of bases of $V$ or $W$. \\
\quad It can be shown that the reduced characters must satisfy the inequality:
\begin{equation}\label{dim_intergral_elem}
 \dim A^{(1)}\leq s_1+2s_2+\dots +ns_n,
\end{equation}
where $A^{(1)}$ is the {\it first prolongation of A}, namely
\begin{equation*}
 A^{(1)}:=(A\otimes V^*)\cap (W\otimes S^2 V^*),
\end{equation*}
and $S^2 V^*$ is the space of symmetric 2-tensors of $V^*$.\\
\quad We reach in this way to one of the most important notion in the theory of exterior differential systems. The tableau   
$A\in  W\otimes V^*$ is called {\it involutive} if equality holds in (\ref{dim_intergral_elem}), i.e.
we have
\begin{equation*}
 \dim A^{(1)}= s_1+2s_2+\dots +ns_n.
\end{equation*}
This condition is also called {\it Cartan test for involutivity}.\\
\quad If $A$ is involutive such that $s_l\neq 0$ and $s_{l+1}= 0$, then $s_l$ is called the {\it character} of the system and 
the integer $l$ is called the {\it Cartan integer} of the system.\\
\quad We can give now the main tool used in this paper, the Cartan--K\"ahler Theorem for Linear Pfaffian systems. Even 
though the theorem can be formulated in general for arbitrary exterior differential systems (see \cite{Br et al 1991}, 
\cite{IL2003}), the version for Linear Pfaffian systems will suffice for our purposes in the present paper.

{\bf Theorem A.1.\ \ The Cartan--K\"ahler Theorem for Linear Pfaffian systems}

\quad {\it Let (I,J) be an analytic linear Pfaffian system on a manifold $\Sigma$, let $x\in \Sigma$ be a point and let $U
\subset \Sigma$ be a neighborhood containing $x$, such that for all $y\in U$,
\begin{enumerate}
\item The apparent torsion is absorbable at $y$, and
\item  the tableau $A_y$ is involutive.
\end{enumerate}
\quad Then solving a series of well-posed Cauchy problems  yields analytic integral manifolds of $(I,J)$ passing through 
x. 
}\\
\quad Informally, one says that the solutions depend (in Cartan-K\"ahler sense) on $s_l$ functions of $l$ variables, 
where $s_l$ is the character of the system (see \cite{IL2003}, p. 176 for the precise statement of the Theorem and other details).
\quad A linear Pfaffian system satisfying the conditions (1), and (2) in the Cartan--K\"ahler Theorem for linear Pfaffian 
systems is said to be {\it involutive}.\\
\quad Recall that if an EDS is not a linear Pfaffian system, then by prolongation one can linearize it and then study its involutivity by 
 Cartan--K\"ahler Theorem for linear Pfaffian systems.

\bigskip

\medskip

\begin{center}

Sorin V. SABAU\\

School of Science, Department of Mathematics\\
Tokai University,\\
Sapporo, 
005\,--\,8601 Japan

\medskip
{\tt sorin@tspirit.tokai-u.jp}

\bigskip

Kazuhiro SHIBUYA\\
Graduate School of Science, Hiroshima University, \\
Higashi Hiroshima, 739\,--\,8521, Japan

\medskip
{\tt shibuya@hiroshima-u.ac.jp}

\bigskip
Hideo SHIMADA\\
School of Science, Department of Mathematics\\
Tokai University,\\
Sapporo, 
005\,--\,8601 Japan

\medskip
{\tt shimadah@tokai-u.jp}

\end{center}

\end{document}